\documentclass[final]{siamltex}
\usepackage{amsfonts,graphicx}
\usepackage[leqno]{amsmath}

\def\jacx{P_n^{(\alpha,\beta)}(x)}
\def\ja1x{P_{n-1}^{(\alpha,\beta)}(x)}
\def\jacw8{W_{\alpha,\beta}}
\def\SI{\mathcal{I}}
\def\eg{{\it e.g.}\  }
\def\ie{{\it i.e.}\  }
\def\f12{\frac{1}{2}}

\newcommand{\DS}{\displaystyle}
\newcommand{\boxit}[1]{\vbox{\hrule\hbox{\strut \vrule \vbox{ \kern6pt #1} \vrule}\hrule}}
\newtheorem{rmk}{Remark}
\title{Spectrum of the {J}acobi tau approximation for the second derivative operator}% with Dirichlet boundary conditions}

\author{Marios Charalambides\thanks{Department of Business Administration, Frederick Institute
of Technology, 7 Yianni Frederickou Street, Pallouriotissa, PO Box
24729, 1303 Nicosia, Cyprus (\texttt{bus.chm@fit.ac.cy}).}
        \and Fabian Waleffe\thanks{Department of Mathematics,
University of Wisconsin, Madison, WI 53706, USA (\texttt{waleffe@math.wisc.edu}). This work was supported in part by NSF grant DMS-0204636.\hfill [Preprint Jul 24, 2006]}}

\date{2006/07/24  FW}  % MARIOS: CHANGE DATE HERE WHEN YOU MAKE CHANGES

\begin{document}

\maketitle

\begin{abstract}
It is proved that the eigenvalues of the Jacobi Tau method for the
second derivative operator with Dirichlet boundary conditions are
real, negative and distinct for a range of the Jacobi parameters.
Special emphasis is placed on the symmetric case of the Gegenbauer
Tau method where the range of parameters included in the theorems
can be extended and characteristic polynomials given by successive
order approximations interlace. This includes the common Chebyshev
and Legendre, Tau and Galerkin methods. The characteristic polynomials for the Gegenbauer Tau method are shown to obey three term recurrences plus a constant term which vanishes for the Legendre Tau and Galerkin cases. These recurrences are equivalent to a tridiagonal plus one row matrix structure. The spectral integration formulation of the Gegenbauer Tau method is shown to lead directly to that fundamental and well-conditioned tridiagonal plus one row matrix structure. A Matlab code is provided.
\end{abstract}

\begin{keywords}
Jacobi polynomials, Gegenbauer polynomials, stable polynomials,
positive pairs, zeros of polynomials, spectral methods
\end{keywords}

\begin{AMS}
65D30, 65L10, 65L15, 65M70, 65N35, 26C10
\end{AMS}

\pagestyle{myheadings}
\thispagestyle{plain}
\markboth{MARIOS CHARALAMBIDES AND FABIAN WALEFFE}{SPECTRUM OF 
JACOBI TAU  SECOND DERIVATIVE OPERATOR}

\section{Introduction}

Constructing polynomial approximations to solutions of differential
equations is the basic ingredient of most numerical methods.
Approximations based on orthogonal polynomials have been
widely used  (\eg \cite{Boyd}, \cite{CHQZ}, \cite{GO77}) because their rate of convergence is faster than
algebraic for arbitrary boundary conditions when the solution is smooth. The purpose of
this paper is to give a rigorous proof that the spectrum of the
Jacobi Tau approximation is real, negative and distinct for the second order
operator with Dirichlet boundary conditions. The Jacobi Tau class of \textit{spectral methods} includes the common Chebyshev and Legendre Tau and Galerkin formulations, as demonstrated below. 
The general method of proof %in section \ref{zeros}  
is similar to that used by Gottlieb and Lustman \cite{G81,GL83} to prove such results for the Chebyshev collocation operator. However, we argue  in section \ref{stabHB} 
that Gottlieb and Lustman's proof for the collocation operator is not complete.

The spectrum of Jacobi Tau approximation for the 1st order operator has been
considered elsewhere \cite{CCW}. Here, we consider polynomial
approximations to the eigenvalue problem
\begin{equation}\label{SOP}
 \frac{d^2u}{dx^2}=\lambda u  \qquad -1<x<1 \qquad {\rm with}
\qquad u(\pm 1)=0.
\end{equation}
The spectrum of Jacobi polynomial approximations to this
eigenvalue problem is directly relevant to numerical simulations
of the diffusion equation $u_t=u_{xx}$ which is itself a building
block for numerical solution of various other problems including
the Stokes and Navier-Stokes equations (\eg \cite[\S 5.1]{CHQZ}, \cite{WT88}).

The 2nd order problem (\ref{SOP}) is a self-adjoint, negative
definite Sturm-Liouville differential eigenproblem, so its
eigenvalues $\lambda$ are real, negative and distinct. The
eigenmodes separate into even and odd modes and have the simple
exact expressions
\begin{equation} %06/7/6
\begin{array}{lll}
 u_e(x) = \cos(2k-1) \frac{\pi}{2} x,    & \DS \lambda = -(2k-1)^2 \frac{\pi^2}{4},
 \cr
\cr
u_o(x)= \sin k \pi x,    & \lambda = -k^2 \pi^2.
\cr
\end{array}
\end{equation}
for $k=1,2,3,\ldots$.

If $u_n(x)$ is a polynomial approximation of degree $n$ to the
exact solution $u(x)$, then $u_n(x)$ satisfies the following
differential equation
\begin{equation}
\lambda u_n(x)-D^2u_n(x)=R_n(x) \label{Rn}
\end{equation}
where the \textit{residual} $R_n(x)$ is a polynomial of degree $n$
in $x$ and $D= d/dx$. We can invert this relation to express the
polynomial approximation $u_n(x)$ in terms of the residual
$R_n(x)$ \cite{GL83}
\begin{equation}\label{u_n}
u_n(x)=\mu \sum_{k=0}^{[n/2]} \mu^k D^{2k} R_n(x)
\end{equation}
where $\mu = 1/\lambda$, and $[n/2]$ denotes the greatest integer
less or equal to $n/2$. We can assume that $\lambda \ne 0$ because
$\lambda=0$ with $u_n(\pm 1)=0$ necessarily corresponds to the
trivial solution $u_n(x)=0$, $\forall x$ in $[-1,1]$, as shown below. 
The inversion (\ref{u_n}) follows from formal
application of the geometric (Neumann) series for 
$(1-\mu D^2)^{-1}= \sum_{k=0}^{\infty} \mu^k D^{2k}$ which terminates
since $R_n(x)$ is a polynomial.
That inversion can also be derived by repeated application of the operator $\mu D^2$ to equation (\ref{Rn}).
Summation of the resulting suite of equations leads to (\ref{u_n}) thanks to telescopic cancelations on the left hand side.
%This inversion is related to the Liouville-Neumann series in more general contexts.

Spectral methods fit in the general framework of the \textit{method
of weighted residuals} \cite{CHQZ}. In the \textit{Tau} method
\cite[\S 10.4.2]{CHQZ}, the polynomial approximation $u_n(x)$ is
determined from the boundary conditions $u_n(\pm 1) = 0$ and the
requirement that $R_n(x)$ is orthogonal to all polynomials
$p_{n-2}(x)$ of degree $n - 2$ or less with respect to a weight
function $W(x)\ge 0$ in the interval $(-1, 1)$
\begin{equation}
\int_{-1}^1 R_n(x) p_{n-2}(x) W(x) dx = 0. \label{jactau}
\end{equation}
%(method of weighted residual).
 These requirements provide $n+1$ equations for the $n+1$ undetermined constants in the
polynomial approximation $u_n(x)$. For the Jacobi weight function
$W_{\alpha,\beta}(x) = (1-x)^{\alpha}(1+x)^{\beta}$, the residual (\ref{Rn}) 
can be written as
\begin{equation}\label{JTR}
R_n(x)=\tau_0 \lambda \, \jacx +\tau_1 \lambda \, \ja1x
\end{equation}
for some $x$-independent coefficients $\tau_0$ and $\tau_1$, where
$\jacx $ is the Jacobi polynomial of degree $n$ (sect.\ \ref{Jacpolys}). This follows from
orthogonality of the Jacobi polynomials in $-1 < x < 1$ with
respect to the Jacobi weight $\jacw8(x)$ which implies
orthogonality of the Jacobi polynomial of degree $k$ to
\textit{any} polynomial of degree $k-1$ or less with respect to
that weight function. Jacobi polynomials are the most general
class of polynomial solutions of a Sturm-Liouville eigenproblem
that is singular at $\pm 1$ as required for faster than algebraic
convergence \cite[\S 9.2.2, \S 9.6.1]{CHQZ}. It is now easy to verify  from (\ref{Rn}) and (\ref{JTR})
that if $\lambda = 0$ then 
$D^2 u_n(x)=R_n(x)=0$ for all $x$ in $(-1,1)$ but the boundary conditions $u_n(\pm 1)=0$  would then
require that $u_n(x)=0$ for all $x$ in $[-1,1]$. Therefore we can assume that $\lambda \ne 0$.

In the \textit{Galerkin} approach, $u_n(x)$ is determined from the
boundary conditions $u_n(\pm 1)=0$ and orthogonality of the
residual $R_n(x)$ to all polynomials of degree $n$ that vanish at
$x=\pm 1$, with respect to a weight function $W(x)\ge 0$. In other
words, the \textit{test functions} are in the same space
(polynomials of degree $n$) as the trial functions and they
satisfy the same boundary conditions. Such polynomials can be
written in the form $(1-x^2) p_{n-2}(x)$ where $p_{n-2}(x)$ is an
arbitrary polynomial of degree $n-2$, and the Galerkin equations can be written as
\begin{equation}
\int_{-1}^1 R_n(x) (1-x^2) p_{n-2}(x) W(x) dx = 0. \label{jacgal}
\end{equation}
 For the Jacobi weight  $W(x)=\jacw8(x)=(1-x)^{\alpha}(1+x)^{\beta}$,
the Galerkin method  is therefore equivalent to the Tau method for
the weight $W_{\alpha+1,\beta+1}(x)$ and the residual controlled by (\ref{jacgal}) has the form
\begin{equation}\label{JGR}
R_n(x)=\tau_0 \lambda \, P_n^{(\alpha+1,\beta+1)}(x) +\tau_1 \lambda \, P_{n-1}^{(\alpha+1,\beta+1)}(x).
\end{equation}
This residual can be written in terms of the derivatives of  $P_{n+1}^{(\alpha+1,\beta+1)}(x)$ and $\jacx$ by making use of (\ref{jacder}). 
  Since we consider a range of
parameters $\alpha$ and $\beta$, the Jacobi-Tau method also includes
some Jacobi-Galerkin methods. 

In the \textit{collocation} approach, $u_n(x)$ is determined from the boundary conditions $u_n(\pm 1)=0$ and enforcing $R_n(x_j)=0$ at the $n-1$ interior Gauss-Lobatto points $x_j$ such that $D P_n^{(\alpha,\beta)}(x_j)=0 $, $j=1,\ldots,n-1$  \cite[\S 2.2]{CHQZ}. The residual (\ref{Rn}) takes the form
\cite[eqn.\ (4.5)]{GL83}
\begin{equation}\label{JCR}
R_n(x)= (A + B x ) \, D \jacx ,
\end{equation}
for some $A$ and $B$ independent of $x$. The collocation residual (\ref{JCR}) is provided for completeness since we do not have results about the collocation method and we raise doubts about the validity of the proof proposed in \cite{GL83}. That residual can be written in several equivalent forms by using the properties of Jacobi polynomials (sect.\ \ref{Jacpolys}).

The characteristic polynomials for the eigenvalues $\mu=1/\lambda$ are derived in section \ref{charpolys} from the explicit expression (\ref{u_n}) for $u_n(x)$ in terms of the residual $R_n(x)$ whose form is specified by the Jacobi Tau (or Galerkin) method as in (\ref{JTR}) and (\ref{JGR}). The zeros of these characteristic  polynomials are shown to be real, negative and distinct 
in section \ref{zeros}.
   Recurrence relations for the Gegenbauer Tau characteristic polynomials are derived in section \ref{rec} where it is shown that the underlying  fundamental matrix structure is tridiagonal + one row. Section \ref{numimp} discusses some implementation issues and shows that the spectral integration implementation directly leads to the tridiagonal + one row structure which is well-conditioned. 
Some of the key properties of Jacobi and Gegenbauer polynomials used in this paper are summarized in appendix \ref{orthopolys}. We use a non-standard normalization for Gegenbauer polynomials, denoted $G_n^{(\gamma)}(x)$, since the standard normalization $C_n^{(\gamma)}(x)$ is singular in the Chebyshev case.

\section{Characteristic Polynomials}
\label{charpolys}

\subsection{Jacobi-Tau method} \label{subjac}

Substituting (\ref{JTR}) into (\ref{u_n}), the Jacobi-Tau approximation can be written
explicitly in terms of the yet undertermined  constants $\tau_0$,
$\tau_1$ and the eigenvalue $\mu=1/\lambda$, as
\begin{equation} \label{Pu_n}
u_n(x)=
\tau_0 \sum_{k=0}^{[\frac{n}{2}]} \mu^{k} D^{2k} \jacx  +
\tau_1 \sum_{k=0}^{[\frac{n-1}{2}]} \mu^{k} D^{2k} \ja1x .
\end{equation}
The boundary conditions $u_n(\pm 1)=0$ then yield the characteristic equations
\begin{equation} \left\{ \begin{aligned}
\tau_0 \sum_{k=0}^{[\frac{n}{2}]} \mu^k D^{2k} P_n^{(\alpha,\beta)}(1) &+
\tau_1 \sum_{k=0}^{[\frac{n-1}{2}]}\mu^k D^{2k} P_{n-1}^{(\alpha,\beta)}(1) &=0, \\
\tau_0 \sum_{k=0}^{[\frac{n}{2}]} \mu^k D^{2k} P_n^{(\alpha,\beta)}(-1) &+
\tau_1 \sum_{k=0}^{[\frac{n-1}{2}]} \mu^k D^{2k} P_{n-1}^{(\alpha,\beta)}(-1) &=0.
\end{aligned} \right.
\end{equation}
Equation (\ref{jacobi}) shows that $P_n^{(\alpha,\beta)}(-1)=(-1)^n P_n^{(\beta,\alpha)}(1)$ 
so the 2nd equation above can be rewritten at $x=1$ by flipping the indices $\alpha$ and $\beta$,
\begin{equation} \left\{ \begin{aligned}
\tau_0 \sum_{k=0}^{[\frac{n}{2}]} \mu^k D^{2k} P_n^{(\alpha,\beta)}(1) &+ 
\tau_1 \sum_{k=0}^{[\frac{n-1}{2}]} \mu^k D^{2k} P_{n-1}^{(\alpha,\beta)}(1) & =0, \\
\tau_0 \sum_{k=0}^{[\frac{n}{2}]} \mu^k D^{2k} P_n^{(\beta,\alpha)}(1) &- 
\tau_1 \sum_{k=0}^{[\frac{n-1}{2}]} \mu^k D^{2k}  P_{n-1}^{(\beta,\alpha)}(1) &=0.
\end{aligned} \right.
\end{equation}
 This system has a non-trivial solution $(\tau_0,\tau_1) \ne(0,0)$ if and only if
\begin{multline}\label{jacpoly}
 %\left(
 \sum_{k=0}^{[\frac{n}{2}]} \mu^k D^{2k}
P_n^{(\alpha,\beta)}(1) \sum_{k=0}^{[\frac{n-1}{2}]} \mu^k D^{2k}
P_{n-1}^{(\beta,\alpha)}(1) \\
+ \, \sum_{k=0}^{[\frac{n-1}{2}]} \mu^k
D^{2k} P_{n-1}^{(\alpha,\beta)}(1)\sum_{k=0}^{[\frac{n}{2}]} \mu^k
D^{2k}
P_n^{(\beta,\alpha)}(1) %\right)
=0.
\end{multline}
This is the characteristic equation for the eigenvalue $\mu$.

\subsection{Gegenbauer-Tau method}

Gegenbauer polynomials $G_n^{(\gamma)}(x)$  are the class of Jacobi polynomials $\jacx$ with equal indices $\alpha=\beta=\gamma-1/2$ (sect.\ \ref{Gegpolys}). The Gegenbauer polynomials are even in $x$ for $n$ even and odd for $n$ odd \cite[22.4.2]{AS}.  Chebyshev and Legendre polynomials are Gegenbauer polynomials with $\gamma=0$ and $1/2$, respectively.
The symmetry of the differential equation (\ref{SOP}) and of the Gegenbauer polynomials allows decoupling of the discrete problem into even and odd solutions. This parity reduction leads to simpler residuals and simpler forms for the corresponding characteristic polynomials. The residual in the
parity-separated Gegenbauer case contains only one term
\begin{equation}
R_n(x)=\tau_0 \lambda \, G_n^{(\gamma)}(x),
\label{RnGeg}
\end{equation}
where $G_n^{(\gamma)}(x)$ is the $n^{th}$ Gegenbauer polynomial and $n$ is even for even solutions and odd for odd solutions.
Substituting (\ref{RnGeg}) into (\ref{u_n}) provides the Gegenbauer-Tau approximation to (\ref{SOP}) in terms of an undertermined constant $\tau_0$ and the eigenvalue $\mu=1/\lambda$
\begin{equation}
u_n(x)=\tau_0 \sum_{k=0}^{[\frac{n}{2}]} \mu^{k} D^{2k}G_n^{(\gamma)}(x).
\end{equation}
The boundary condition $u_n(\pm 1)=0$ leads to the characteristic
polynomial equation
\begin{equation}
\sum_{k=0}^{[\frac{n}{2}]} \mu^k D^{2k} G_n^{(\gamma)}(1)=0,
\label{gegchara}
\end{equation}
since by symmetry $G_n^{(\gamma)}(-1)=(-1)^n G_n^{(\gamma)}(1)$ and
the two boundary conditions give the same equation.

\section{Zeros of characteristic polynomials}
\label{zeros}

\subsection{Stable polynomials and the Hermite Biehler Theorem}

\label{stabHB}

The general approach to prove that the eigenvalues are real, negative and distinct  is to construct a particular \textit{stable} polynomial $p(z)$ then to use the Hermite Biehler theorem to deduce that the polynomials $\Omega_1(\mu)$ and $\Omega_2(\mu)$ such that $p(z) = \Omega_1(z^2) + z \Omega_2(z^2)$ have real, negative and distinct zeros that interlace. 

\begin{definition} \label{defstable}
A real polynomial, $p(z)$, is  a stable polynomial
(or a {\it Hurwitz polynomial}), if all its zeros  lie in the open
left half-plane, \ie their real part is strictly less than zero,
$\Re z < 0$. 
\end{definition}
\begin{definition}\label{def1}
Let $\Omega_1(\mu)$ and $\Omega_2(\mu)$ be two real polynomials of
degree $n$ and $n-1$ (or $n$) respectively, then
$\Omega_1(\mu)$ and $\Omega_2(\mu)$ form a positive pair if:
(a)  the roots $\mu_1,\ldots,\mu_{n}   $ of $\Omega_1$ and
$\mu'_1\cdots,\mu'_{n-1}$ (or
$\mu'_1,\cdots, \mu'_{n}  $) of $\Omega_2$ are real, negative and distinct; 
(b) the roots strictly interlace (or alternate) as follows:
$$\mu_1<\mu'_1<\cdots<\mu'_{n-1}<\mu_{n}   <0 \quad  ( \mbox{or}\; 
\mu'_1<\mu_1<\cdots <\mu'_{n}<\mu_{n}   <0 ); $$
(c) the highest coefficients of $\Omega_1(\mu)$ and
$\Omega_2(\mu)$ are of like sign.
\end{definition}

We will use the following theorems about positive pairs, \cite[p.\ 198]{RS}, \cite[Sec.\ 2]{WA}:
\begin{lemma}
\label{PPLClemma} Any nontrivial real linear combination of two
polynomials that form a positive pair has real roots.
\end{lemma}
\begin{lemma} 
\label{WagnerLemma} Let $P(\mu)$ and $Q(\mu)$ be real standard polynomials (\ie the
leading coefficient is positive) %in $R[\mu]$ 
with only non-positive zeros. Then $P(\mu)$ interlaces (or alternates)
$Q(\mu)$ (in the sense of definition \ref{def1}, but not
strictly) if and only if for all $A>0$ both $ Q(\mu)+AP(\mu)$ and
$Q(\mu) + A \mu P(\mu)$ have only non positive zeros.
\end{lemma}

%\begin{rmk}\label{rmk1}
The polynomials $P$ and $Q$ of theorem \ref{WagnerLemma} do not necessarily
form a \textit{positive pair} since they are allowed to have
common and/or multiple roots. We call this set of polynomials a
\textit{quasi-positive pair.}
%\end{rmk}
%Another useful theorem \cite{GL83} is the following:
\begin{lemma} {\rm \cite[Lemma 3.4]{GL83}}
\label{GLlemma} If $\Omega_1(\mu),\Omega_2(\mu)$ and
$\Theta_1(\mu),\Theta_2(\mu)$ are two positive pairs then the zeros of 
%\begin{equation}
$H(\mu)=\Omega_1(\mu)\Theta_2(\mu)+\Omega_2(\mu)\Theta_1(\mu)$
%\end{equation}
are real, negative and distinct.
\end{lemma}

Stability (definition \ref{defstable}) is very important in temporal discretizations
and matrix theory \cite{CHQZ} as well as in analysis (\eg \cite{holtz} and references therein). Stable polynomials can surface as
characteristic polynomials of a numerical method applied on a
differential equation. A necessary and sufficient condition for a
polynomial to be stable is given by the Routh-Hurwitz theorem (see,
for example, \cite[\S 40]{marden},  \cite[\S 23]{O63}). Other
important characterizations of stable polynomials are the
Routh-Hurwitz criterion and the total positivity of a Hurwitz
matrix \cite{holtz}, although these will not be used here.
The characterization of stable polynomials that will be most useful here is given by
the \emph{Hermite-Biehler Theorem} \cite[p.\ 197]{RS},\cite{holtz}.

\begin{theorem}[Hermite-Biehler]\label{HBthm}
The polynomial with real coefficients %\begin{equation}
$p(z)\,=\,\Omega_1(z^2)+z \Omega_2(z^2)$ %\end{equation}
is stable, if and only if $\Omega_1(\mu)$ and $\Omega_2(\mu)$ form
a positive pair.
\end{theorem}

The Hermite Biehler theorem states that 
the even and odd parts of stable polynomials form positive pairs.
This supplies us with a very strong tool to prove reality 
and negativity of the roots of certain polynomials. 

Gottlieb and Lustman \cite{GL83} used the Hermite Biehler theorem to prove that the spectrum of the Chebyshev collocation operator for the heat equation is real, negative and distinct for a variety of homogeneous boundary conditions.  The basic strategy is to show that the characteristic polynomial for that method are the even or odd parts of a stable polynomial. 
  Our results extend their strategy to a class of Jacobi and Gegenbauer Tau methods that includes Chebyshev and Legendre Tau and Galerkin formulations.  Although the general approach is similar to that of Gottlieb \cite{G81} and Gottlieb and Lustman \cite{GL83}, the extension is technically non-trivial and there are differences and some corrections. The key steps in \cite{GL83} is to prove that the polynomials \cite[(4.11),(4.12)]{GL83} are stable. To do so, Gottlieb and Lustman  derive a first order differential equation for those polynomials then transform that ODE into an inhomogeneous one-way wave equation \cite[(4.13)]{GL83} and call on the results \cite[(3.18),(3.20)]{G81} to deduce stability. This is not quite correct since the eigenvalue $\mu$ here is complex hence $w_N(x,t)$ in \cite[(4.13)]{GL83} is also complex while Gottlieb implicitly assumes reality of $v_N(x,t)$ and $R_N(x,t)$ in \cite[(3.18),(3.20)]{G81}. 
  
  The proof for  \cite[(4.11)]{GL83} can be fixed and generalized as done in  \cite{CCW} (theorem \ref{CCWthm} below) where we deduce stability of the polynomials (\ref{Phin}) below without going back to a one-way wave equation. Here, that proof  is further  generalized  (appendix \ref{2proofs}) to the polynomials (\ref{Phi1A}) and (\ref{Phi2A}) below in order to prove our results about the Jacobi Tau method for the 2nd order operator. Our proof follows Gottlieb's ideas to derive the results \cite[(3.18), (3.20)]{G81} although we do not use Gauss integration. 
  
 The proof for \cite[(4.12)]{GL83} does not appear to be correct however and we have not succeeded in obtaining a corrected proof. Gottlieb and Lustman do not provide a proof of stability for that polynomial (4.12), they state only that a \textit{`similar argument holds'}. Gottlieb \cite{G81} likewise suggests that the proof of stability for \cite[(3.11b)]{G81} implies stability for \cite[(3.11a)]{G81} but this is not evident since $\tau_1(t)$ and $\tau_2(t)$ are distinct functions of time that are fully determined by their respective solution procedure.  Gottlieb also suggests that  $v_N(x,t)$ \cite[(3.8)]{G81} is directly related to $u_N(x,t)$ \cite[(3.2)]{G81}  by relation \cite[(3.8)]{G81}. It is true that $T_N(x_n)$ can be eliminated for $n=0,\ldots, N-1$ as suggested in the derivation of \cite[(3.8)]{G81}, since $(1-x) DT_N(x) = 2(-1)^{N-1} \sum_{k=0}^N (-1)^k c_k^{-1} T_k(x)$ as can be deduced from \cite[(3.2),(3.3)]{G81}. However this does not imply that the resulting $d_k$ coefficients \cite[(3.8)]{G81} deduced from the $a_k$'s  that solve \cite[(3.6)]{G81} are the same $d_k$'s as those that  solve  \cite[(3.10)]{G81}.
 
 Hence, it appears that there is currently no proof of the stability of \cite[(3.6),(3.11a)]{G81} and \cite[(4.12)]{GL83}, therefore invalidating Gottlieb and Lustman's proof that the eigenvalues of the Chebyshev collocation operator are real, negative and distinct \cite{GL83}.
  
% Parity-reduced Gegenbauer collocation gives
% $R_n= x D G_n^{(\gamma)} (x)$. From (\ref{gegder}), (\ref{gegrec}) and (\ref{ggp}) one deduces that 
% $x D G_n^{(\gamma)} (x) = 2 (\gamma+1) G_n^{(\gamma+1)} (x) -(2\gamma+n) G_n^{(\gamma)} (x)$ hence the Gegenbauer collocation residual is one particular linear combination of the Galerkin $G_n^{(\gamma+1)} (x)$ and Tau $G_n^{(\gamma)} (x)$ residuals. 

\subsection{Important Stable Polynomials and Positive Pairs}

Here we prove stability of certain real polynomials whose even and odd parts
are directly related to the characteristic polynomials derived in section \ref{charpolys} for the Jacobi Tau method.
\begin{theorem}  \label{CCWthm}
 Let $P_n^{(\alpha,\beta)}(x)$ denote the Jacobi
polynomial of degree $n$, where $n\ge 2$. If $-1< \alpha \le1$ and
$\beta>-1$, then the zeros of the polynomial
\begin{equation}
\Phi_n(\mu):= \sum_{k=0}^{n}\left(\frac{d^k}{d
x^k}P_n^{(\alpha,\beta)}(x)\right )_{x=1} \, \mu^k \label{Phin}
\end{equation}
lie in the left half-plane; that is, $\Phi_n(\mu)$ is a stable
polynomial.
\end{theorem}
The proof of this theorem is in \cite{CCW} together with a
discussion of its relation to zeros of Bessel polynomials. 
The next two theorems give two generalizations of the
above result that are needed for this paper.
\begin{theorem}
Let $P_n^{(\alpha,\beta)}(x)$ denote the Jacobi polynomial of degree
$n$, with $n\ge 3$, then
the polynomial \label{thmP01}
\begin{equation}
\Phi^1_n(\mu):=\sum_{k=0}^{n} \left( \frac{d^k}{dx^k}
P_n^{(\alpha,\beta)}(x)\right)_{x=1} \mu^k+A\sum_{k=0}^{n-1}
\left( \frac{d^k}{dx^k} P_{n-1}^{(\alpha,\beta)}(x) \right)_{x=1}
\mu^k
\label{Phi1A}
\end{equation}
is stable for every $A \ge 0$ when $-1< \alpha \le 0$ and $\beta > -1$.
\end{theorem}

\begin{theorem}
Let $P_n^{(\alpha,\beta)}(x)$ denote the Jacobi polynomial of degree
$n$, with $n\ge 3$, then
the polynomial \label{thmP02}
\begin{equation}
\Phi^2_n(\mu):=\sum_{k=0}^{n} \left( \frac{d^k}{dx^k}
P_n^{(\alpha,\beta)}(x)\right)_{x=1} \mu^k+A \mu^2
\sum_{k=0}^{n-1} \left( \frac{d^k}{dx^k}
P_{n-1}^{(\alpha,\beta)}(x) \right)_{x=1} \mu^k
\label{Phi2A}
\end{equation}
is stable for every $A \ge 0$ when  $-1< \alpha \le 1$ and $\beta > -1$.
\end{theorem}

The proofs of these theorems are technical and they are given in
appendix \ref{2proofs}. Our next theorem combines all the above theorems to
get an important result.
\begin{theorem}\label{JacPP}
Let $n \ge 3$, then the polynomials
\begin{equation}
\Omega_n^{(\alpha,\beta)}(\mu):=\sum_{k=0}^{[\frac{n}{2}]} \mu^k
D^{2k} P_n^{(\alpha,\beta)}(1), \quad
\Omega_{n-1}^{(\alpha,\beta)}(\mu):=\sum_{k=0}^{[\frac{n-1}{2}]}
\mu^k D^{2k} P_{n-1}^{(\alpha,\beta)}(1)
\end{equation}
form a positive pair if $-1 <\alpha,\beta \le 0$ or $0 <\alpha,\beta \le 1$.
\end{theorem}

\begin{rmk}
It was shown in \cite{CCW} that these polynomials
have real negative and distinct roots for  $-1 <\alpha\le 1$ and $-1 <\beta $. 
The important addition of theorem \ref{JacPP}  is that the roots of these polynomials
interlace as was conjectured in \cite{CCW}.
\end{rmk}

\begin{proof}
Applying the \textit{Hermite-Biehler} Theorem to the stable
polynomials of theorem \ref{CCWthm} for a given $n$ and also for $n-1\ge 2$, proves that the polynomials
$\Omega_n^{(\alpha,\beta)}(\mu)$ and
$\Omega_{n-1}^{(\alpha,\beta)}(\mu)$ have real negative and distinct
roots for $-1<\alpha ,\beta \le 1$ (Notice that we can interchange
$\alpha$ and $\beta$). Applying the \textit{Hermite-Biehler}
Theorem to the stable polynomials of theorems \ref{thmP01} and
\ref{thmP02} shows that the polynomials
\begin{gather}
\sum_{k=0}^{[\frac{n}{2}]} \mu^k D^{2k}
P_n^{(\alpha,\beta)}(1)+A\sum_{k=0}^{[\frac{n-1}{2}]} \mu^k D^{2k}
P_{n-1}^{(\alpha,\beta)}(1)=\Omega_n^{(\alpha,\beta)}(\mu)+A\Omega_{n-1}^{(\alpha,\beta)}(\mu), \\
\sum_{k=0}^{[\frac{n}{2}]} \mu^k D^{2k} P_n^{(\alpha,\beta)}(1) +A
\mu \sum_{k=0}^{[\frac{n-1}{2}]} \mu^k D^{2k}
P_{n-1}^{(\alpha,\beta)}(1)=\Omega_{n}^{(\alpha,\beta)}(\mu)+A\mu\Omega_{n-1}^{(\alpha,\beta)}(\mu)
\end{gather}
have real negative and distinct roots for all $A>0$ and $-1
<\alpha,\beta \le 0$.  These results provide sufficient information
to apply lemma \ref{WagnerLemma} to deduce  that the set of polynomials
$(\Omega_n^{(\alpha,\beta)}(\mu),\Omega_{n-1}^{(\alpha,\beta)}(\mu))$
form quasi-positive pair for $-1<\alpha ,\beta \le 0$.

To show that these polynomials form a positive pair
 recall that both $\Omega_n^{(\alpha,\beta)}(\mu)$ and
$\Omega_{n-1}^{(\alpha,\beta)}(\mu)$ have real, negative and
distinct roots by theorem \ref{CCWthm}. Thus it remains to show that
they have no common roots. To do so, assume that $\mu$ is a common root. Then
$P_1(\mu)=\Omega_n^{(\alpha,\beta)}(\mu)+A\Omega_{n-1}^{(\alpha,\beta)}(\mu)=0
$ and $P_2(\mu)=\Omega_n^{(\alpha,\beta)}(\mu)+A\mu
\Omega_{n-1}^{(\alpha,\beta)}(\mu)=0 $. Since
$D\Omega_n^{(\alpha,\beta)}(\mu)\not = 0$ and
$D\Omega_{n-1}^{(\alpha,\beta)}(\mu) \not = 0$ (both $\Omega_n$ and
$\Omega_{n-1}$ do not have a double root), set $A=
-\frac{D\Omega_{n}^{(\alpha,\beta)}(\mu)}{D\Omega_{n-1}^{(\alpha,\beta)}(\mu)}$
or $A= -\frac{D\Omega_{n}^{(\alpha,\beta)}(\mu)}{\mu
D\Omega_{n-1}^{(\alpha,\beta)}(\mu)}$, whichever one is positive
(one of the two must be since $\mu$ is negative). But this will
imply that $DP_1(\mu)=0 $ and $DP_2(\mu)=0 $ respectively, a
contradiction since $P_1(\mu)$ and $P_2(\mu)$ have simple zeros. For
the second range of parameters, replace $n$ with $n+1$ in theorems
\ref{thmP01} and \ref{thmP02} and apply the \textit{Hermite-Biehler}
Theorem. This gives that the polynomials
\begin{gather}
\sum_{k=0}^{[\frac{n}{2}]} \mu^k D^{2k+1}
P_{n+1}^{(\alpha,\beta)}(1)+A\sum_{k=0}^{[\frac{n-1}{2}]} \mu^k
D^{2k+1} P_{n}^{(\alpha,\beta)}(1), \\
%\end{equation} and \begin{equation}
\sum_{k=0}^{[\frac{n}{2}]} \mu^k D^{2k+1}
P_{n+1}^{(\alpha,\beta)}(1) +A \mu \sum_{k=0}^{[\frac{n-1}{2}]}
\mu^k D^{2k+1} P_{n}^{(\alpha,\beta)}(1)
\end{gather}
also have real negative and distinct roots. Using (\ref{jacder}) these polynomials transform to
\begin{gather}
\sum_{k=0}^{[\frac{n}{2}]} \mu^k D^{2k}
P_n^{(\alpha+1,\beta+1)}(1)+A'\sum_{k=0}^{[\frac{n-1}{2}]} \mu^k
D^{2k} P_{n-1}^{(\alpha+1,\beta+1)}(1), \\
%\end{equation} and \begin{equation}
\sum_{k=0}^{[\frac{n}{2}]} \mu^k D^{2k} P_n^{(\alpha+1,\beta+1)}(1)
+A' \mu \sum_{k=0}^{[\frac{n-1}{2}]} \mu^k D^{2k}
P_{n-1}^{(\alpha+1,\beta+1)}(1).
\end{gather}
The new polynomials have real negative and distinct roots for all
$A'>0$, where $A'=\frac{n+\alpha+\beta}{n+1+\alpha+\beta} A$, and
thus the second range $0 <\alpha,\beta \le 1$ follows from the first
one by a simple change of variables. %This completes the proof.
\end{proof}

\subsection{Eigenvalues of the Gegenbauer and Jacobi Tau methods}

The previous subsection provides all necessary information needed for
proving reality and negativity of the eigenvalues. First consider
the Gegenbauer case.

\begin{theorem}\label{thmG}
The eigenvalues of the Gegenbauer Tau discretization of the second
order operator with Dirichlet boundary conditions, problem
\ref{SOP}, are real negative and distinct for $-1/2<\gamma \le 5/2$.
Also, characteristic polynomials given by successive order (\ie $n-1$ and $n$)
approximations interlace.
\end{theorem}

\begin{proof}
The eigenvalues are the roots of (\ref{gegchara}). The theorem follows directly
from theorem \ref{JacPP} since $\alpha=\beta=\gamma-1/2$
in the Gegenbauer case and the two ranges for the indices $(\alpha, \beta)$ merge into the single range $-1/2 < \gamma \le 5/2$. 
\end{proof}

\begin{rmk} The results of theorem \ref{thmG} are sharp in the sense that well-conditioned numerical calculations (sect.\ \ref{numerical}) give some complex conjugate pairs of eigenvalues for $\gamma>5/2$ and Gegenbauer integration (\ref{gegortho}) diverges in general for $\gamma \le -1/2$.

\end{rmk} 
\begin{rmk}\label{Classical}
As pointed out in the introduction, the Galerkin method with weight
function $W_{\alpha,\beta}(x)$ for problem \ref{SOP} is equivalent
to the Tau method with weight function $W_{\alpha+1,\beta+1}(x)$. 
For the Gegenbauer case the Galerkin method with weight
function $W_{\gamma}(x)$ is equivalent to the Tau method with weight
function $W_{\gamma+1}(x)$ . Thus, a direct consequence of theorem
\ref{thmG}  is that the eigenvalues of the Gegenbauer Galerkin
method are real negative and distinct for $-1/2<\gamma \le 3/2$.
Again, characteristic polynomials given by successive order
approximations interlace.
\end{rmk}

\begin{rmk}
Theorem \ref{thmG} includes the  Chebyshev and Legendre polynomials. For $\gamma=0$, we have $G_n^{(0)}(x) =\frac{T_n(x)}{n}$, where
$T_n(x)$ denotes the $n^{th}$ Chebyshev polynomial of the first kind
\cite[p.\ 19]{RS}. Thus the theorem implies that the Chebyshev-Tau
method has real negative and distinct eigenvalues with interlacing
characteristic polynomials given by successive order approximations.
For $\gamma=\frac{1}{2}$ we have $G_n^{(1/2)}(x)=P_n(x)$, where
$P_n(x)$ is the $n^{th}$ Legendre polynomial and for $\gamma=1$ we
have that $G_n^{(1)}(x)=\frac{U_n(x)}{2}$ where $U_n(x)$ the
$n^{th}$ Chebyshev polynomial of second kind. Therefore, the same
result holds for both the Legendre Tau and the Chebyshev Tau of 2nd
kind method. Furthermore, in the Galerkin case (see remark
\ref{Classical}) the Galerkin Chebyshev, the Galerkin Legendre and the
Galerkin Chebyshev of the 2nd kind methods have real, negative and distinct eigenvalues as well.
\end{rmk}

For the Jacobi case we have
\begin{theorem}\label{ThmJac}
The eigenvalues of the Jacobi Tau discretization of the second order
operator with Dirichlet boundary conditions, problem \ref{SOP}, are
real negative and distinct if  $-1 <\alpha,\beta \le 0$ or $0 <\alpha,\beta \le 1$.
\end{theorem}

\begin{proof}
By theorem \ref{JacPP} the polynomials
$(\Omega_n^{(\alpha,\beta)}(\mu),\Omega_{n-1}^{(\alpha,\beta)}(\mu))$
form a positive pair for $-1 <\alpha,\beta \le 0$ and $0
<\alpha,\beta \le 1$. Interchanging the indices $\alpha$, $\beta$, the
same result holds for 
$(\Omega_n^{(\beta,\alpha)}(\mu),\Omega_{n-1}^{(\beta,\alpha)}(\mu))$.
Application of theorem \ref{GLlemma} to these sets of positive pairs
gives that the polynomial
\begin{equation}
B_n^{(\alpha,\beta)}(\mu)=\Omega_n^{(\alpha,\beta)}(\mu)\Omega_{n-1}^{(\beta,\alpha)}(\mu)+\Omega_n^{(\beta,\alpha)}(\mu)\Omega_{n-1}^{(\alpha,\beta)}(\mu)
\end{equation}
has real negative and distinct roots. Equation
(\ref{jacpoly}) shows that $B_n^{(\alpha,\beta)}(\mu)$ is the
characteristic polynomial for the Jacobi Tau method.
\end{proof}

This paper focuses on the second order problem with Dirichlet
boundary conditions. Naturally, questions arise about the equivalent results for 
different boundary conditions. The next remarks gives some answer to that question.
\begin{rmk}
Consider problem \ref{SOP} with Neumann boundary conditions \ie
$\lambda u-D^2u=0$ with $Du(\pm 1)=0$. The Jacobi Tau method gives
$\lambda u-D^2u=\tau_0 \jacx +\tau_1 \ja1x$. Notice that $\lambda=0$
is an eigenvalue since $u=constant$ is a solution. For $\lambda \not
= 0$, differentiate the last equation to get $\lambda Du-D^3u=\tau_0
D\jacx +\tau_1 D\ja1x$ with $Du(\pm 1)=0$. Set now $v=Du$ and with
the use of (\ref{jacder}) the equation transforms to $\lambda
v-D^2v=\tau_0' P_{n-1}^{(\alpha+1,\beta+1)}(x) +\tau_1'
P_{n-2}^{(\alpha+1,\beta+1)}(x) $ with $v(\pm 1)=0$. This is the
Jacobi Tau approximation of second kind \ie
$(\alpha,\beta)\rightarrow (\alpha+1,\beta+1)$ for problem
\ref{SOP}, and by theorem \ref{ThmJac} it has real negative and
distinct eigenvalues for $-1 <\alpha,\beta \le 0$.
\end{rmk}
\begin{rmk}
Consider problem \ref{SOP} with boundary conditions $u(-1)=0$ and
$Du(1)=0$. Using equation \ref{Pu_n} and following the ideas of
subsection \ref{subjac} we get that the characteristic polynomial
in this case is $B_n^{(\alpha,\beta)}(\mu)=k_{n-1}
\Omega_n^{(\beta,\alpha)}(\mu)\Omega_{n-2}^{(\alpha+1,\beta+1)}(\mu)+k_{n}\Omega_{n-1}^{(\beta,\alpha)}(\mu)\Omega_{n-1}^{(\alpha+1,\beta+1)}(\mu)$,
where $k_n= \frac{1}{2} (n+\alpha+\beta+1)$. Since the polynomials
$(\Omega_n^{(\beta,\alpha)}(\mu),\Omega_{n-1}^{(\beta,\alpha)})$
and
$(k_{n}\Omega_{n-1}^{(\alpha+1,\beta+1)}(\mu),k_{n-1}\Omega_{n-2}^{(\alpha+1,\beta+1)})$
form positive pairs for $-1 <\alpha,\beta \le 0$ then by theorem
\ref{GLlemma} the roots of $B_n^{(\alpha,\beta)}(\mu)$ are real
negative and distinct for $-1 <\alpha,\beta \le 0$.

An alternative implementation is to rescale the domain to $[-1,0]$ with the Neuman boundary condition $Du(0)=0$ and to use only the even Gegenbauer polynomials with a Gegenbauer-Tau approach on the entire domain $[-1,1]$ and boundary conditions $u(\pm 1)=0$ since such polynomials automatically satisfy $Du(0)=0$.  
\end{rmk}

\section{Characteristic polynomial recurrences}% Relations}
\label{rec}

The previous section shows that the characteristic polynomials given by successive order approximations of the Gegenbauer Tau method have real negative and distinct roots that
interlace strictly, provided that $-1/2 < \gamma \le 5/2$. Reality of the roots as well as the interlacing
property is an important characteristic of orthogonal
polynomials of successive order \cite{szego},\cite[p.\ 16]{RS}. These properties are
direct consequences of %results that are based on
 the three term recurrence relations satisfied by orthogonal polynomials %, such as the Christoffel-Darboux formula 
  \cite{szego},\cite{RS}.
Here, we derive recurrence relations for the characteristic polynomials of the Gegenbauer Tau
method and show that they consist of three term recurrences {\em plus} a constant term in general.
The constant term vanishes when $\gamma=1/2$ or $3/2$.

From (\ref{gegchara}), the characteristic polynomials  for  Gegenbauer-Tau approximations to the even and odd modes of (\ref{SOP}) are, respectively,
\begin{equation}
p_{m}(\mu) := \sum_{k=0}^{m} \mu^k  \,D^{2k} G^{(\gamma)}_{2m}(1),
\quad \mbox{and} \quad
q_{m}(\mu) := \sum_{k=0}^{m} \mu^k  \, D^{2k} G^{(\gamma)}_{2m+1}(1).
\label{pmqmdefs}
\end{equation}
theorem \ref{thmG} states that these polynomials have real, negative and distinct zeros for 
$-1/2 < \gamma \le 5/2$ and that the zeros of $p_m(\mu)$ and $q_m(\mu)$ interlace, as do the zeros of $q_{m-1}(\mu)$ and $p_m(\mu)$. 
Recurrence relations for these characteristic polynomials follow directly from recurrences for 2nd derivatives of Gegenbauer polynomials which, using (\ref{dgegrec}) twice, read
\begin{equation}
\left\{\begin{gathered}
G_{0}^{(\gamma)}= \frac{D^2G_{2}^{(\gamma)}}{2(\gamma+1)}, \qquad
G_{1}^{(\gamma)}=
\frac{D^2G_{3}^{(\gamma)}}{4(\gamma+1)(\gamma+2)},  \hfill
%\; G_{2}^{(\gamma)}=  \frac{D^2G_{4}^{(\gamma)}}{4(\gamma+2)(\gamma+3)}       -\frac{D^2G_{2}^{(\gamma)}}{2(\gamma+1)(\gamma+3)},
 \\ %Checked! 
G_n^{(\gamma)}= \hfill \\
\frac{D^2G_{n+2}^{(\gamma)}}{4(\gamma+n+1)(\gamma+ n)} -
\frac{D^2G_{n}^{(\gamma)}}{2(\gamma+ n+1)(\gamma+n-1)} +
\frac{D^2G_{n-2}^{(\gamma)}}{4(\gamma+n)(\gamma+n-1)}. 
\end{gathered}
\right.
\label{dgegrec2}
\end{equation}

\subsection{Recurrences for even modes}

% and $p_1(\mu)=(2\gamma+1)/2 + 2(1+\gamma) \,\mu$ and
Substituting  (\ref{dgegrec2}) with  $n=2m$ into the characteristic polynomials
$p_m(\mu)$ defined in (\ref{pmqmdefs}),  with $p_0(\mu)=1$ since $G_0^{(\gamma)} :=1$,  yields the recurrence
\begin{equation}
\left\{ \begin{gathered}
\mu p_0(\mu)= \frac{p_1(\mu)}{2(\gamma+1)} - K_0^{(\gamma)},  \hfill \\
\mu p_1(\mu)= \frac{p_2(\mu)}{4(\gamma+2)(\gamma+3)}-\frac{p_1(\mu)}{2(\gamma+1)(\gamma+3)} -
K_2^{(\gamma)}, \hfill \\
\mu p_{m}(\mu) = \hfill \\
 \hfill \frac{ p_{m+1}(\mu)}{4(\gamma\!+\!n\!+\!1)(\gamma\!+\!n)} 
-\frac{ p_{m}(\mu)}{2(\gamma\!+\!n\!+\!1)(\gamma\!+\!n\!-\!1)}
+\frac{p_{m-1}(\mu)}{4(\gamma\!+\!n)(\gamma\!+\!n\!-\!1)}  -K_n^{(\gamma)}
\end{gathered}
\right.
\label{pmrec}
\end{equation}
where, using (\ref{Gat1}), 
\begin{equation}
 \begin{aligned}
K_0^{(\gamma)} &:= \frac{G_{2}^{(\gamma)}(1)}{2(\gamma+1)}
      =\frac{2\gamma+1}{4 (\gamma+1)}, \\  %\quad \mbox{and} \quad
K_2^{(\gamma)} &:=\frac{G_{4}^{(\gamma)}(1)}{4(\gamma+2)(\gamma+3)}
      -\frac{G_{2}^{(\gamma)}(1)}{2(\gamma+1)(\gamma+3)}
      =\frac{(2\gamma^2+\gamma-7)(1+2\gamma)}{48(\gamma+1)(\gamma+2)}.
      \label{K0K2}
\end{aligned}
\end{equation}
and
\begin{equation}
K_n^{(\gamma)} =
\frac{G_{n+2}^{(\gamma)}(1)}{4(\gamma+n+1)(\gamma+ n)}
-\frac{G_{n}^{(\gamma)}(1)}{2(\gamma+n+1)(\gamma+n-1)}
+\frac{G_{n-2}^{(\gamma)}(1)}{4(\gamma+n)(\gamma+n-1)}.
\label{Kna}
\end{equation}
for $n\ge 3$. From (\ref{Gat1}), the latter expression reduces to
 \begin{equation}\label{Kn}
K_n^{(\gamma)}=  \frac{(2\gamma-1)(2\gamma-3)}{n(n^2-1)(n+2)} \,
G_{n-2}^{(\gamma)}(1)
=\frac{(2\gamma-1)(2\gamma-3)}{n(n^2-1)(n^2-4)} \binom{2\gamma+n-3}{n-3},
\end{equation}
 for $n \ge 3$, and obeys the recurrence
\begin{equation}
K_{n+2}^{(\gamma)}=
\frac{(2\gamma+n-1)(2\gamma+n-2)}{(n+4)(n+3)}K_n^{(\gamma)}  \quad
\mbox{with} \quad
K_4^{(\gamma)}=\frac{(4\gamma^2-1)(2\gamma-3)}{720}. \label{Knrec}
\end{equation}

\subsection{Recurrences for odd modes}

Substituting  (\ref{dgegrec2}) with  $n=2m+1$ into the characteristic polynomials for the odd modes $q_{m}(\mu)$ defined in (\ref{pmqmdefs}), with  $q_0(\mu)=1$, gives
%  and $q_1(\mu)=4(\gamma+1)(\gamma+2)\mu + (\gamma+1)(2\gamma+3)$. Using
\begin{equation}
\left\{ \begin{gathered}
\mu q_0(\mu)= \frac{q_1(\mu)}{4(\gamma+1)(\gamma+2)} -
K_1^{(\gamma)} ,  \hfill \\
\mu q_{m}(\mu) = \hfill \\ \frac{ q_{m+1}(\mu)}{4(\gamma\!+\!n\!+\!1)(\gamma\!+\!n)} 
-\frac{ q_{m}(\mu)}{2(\gamma\!+\!n\!+\!1)(\gamma\!+\!n\!-\!1)}
+\frac{q_{m-1}(\mu)}{4(\gamma\!+\!n)(\gamma\!+\!n\!-\!1)}  -K_n^{(\gamma)}
\end{gathered}
\right.
\label{qmrec}
\end{equation}
where, using (\ref{Gat1}), 
$K_1^{(\gamma)} :=[4(\gamma+1)(\gamma+2)]^{-1} G_3^{(\gamma)}(1) = [12 (\gamma+2)]^{-1} (2\gamma+1)$ and  $K_n^{(\gamma)}$ as in (\ref{Kna}) but here with $n=2m+1$. The recurrence (\ref{Knrec}) applies here also but starting now with $K_3^{(\gamma)} = (2\gamma-1)(2\gamma-3)/120$.

In general, (\ref{pmrec}) and (\ref{qmrec}) are three term recurrences \emph{plus} the constants $K_n^{(\gamma)}$. These constants vanish for all $n\ge 3$ when $\gamma=1/2$, the Tau-Legendre method, and when $\gamma=3/2$, the Tau-Legendre of the 2nd kind or Galerkin-Legendre method. In those cases, the recurrences have only three terms hence the corresponding $p_m(\mu)$ and $q_m(\mu)$ sequences of polynomials are orthogonal polynomials \cite[p.\ 13 and references therein]{RS}.   The  recurrence (\ref{Knrec}) for $K_n^{(\gamma)}$ also indicates why $\gamma=5/2$ is a critical value in theorem \ref{thmG}. For $\gamma < 5/2$, $(2\gamma+n-1)(2\gamma+n-2) < (n+4)(n+3)$ and $K_n^{(\gamma)}$ decreases with increasing $n$, while for $\gamma >5/2$, $K_n^{(\gamma)}$ increases with $n$. If $K_n^{(\gamma)}=0$ for $n \ge 3$, the characteristic polynomial sequences $p_m(\mu)$ and $q_m(\mu)$ satisfy a three term recurrence, respectively, therefore they are orthogonal and have real roots that interlace. The constants $K_n^{(\gamma)} \ne 0$  pulls down or pushes up the successive polynomials in the sequences with respect to that orthogonal case. For $\gamma > 5/2$ that shift leads to the bifurcation from real eigenvalues to complex conjugate pairs. 

\section{Numerical Implementation}
\label{numimp}

\subsection{Matrix formulation of the recurrences}
\begin{figure}[h]
\includegraphics[width=0.5\textwidth]{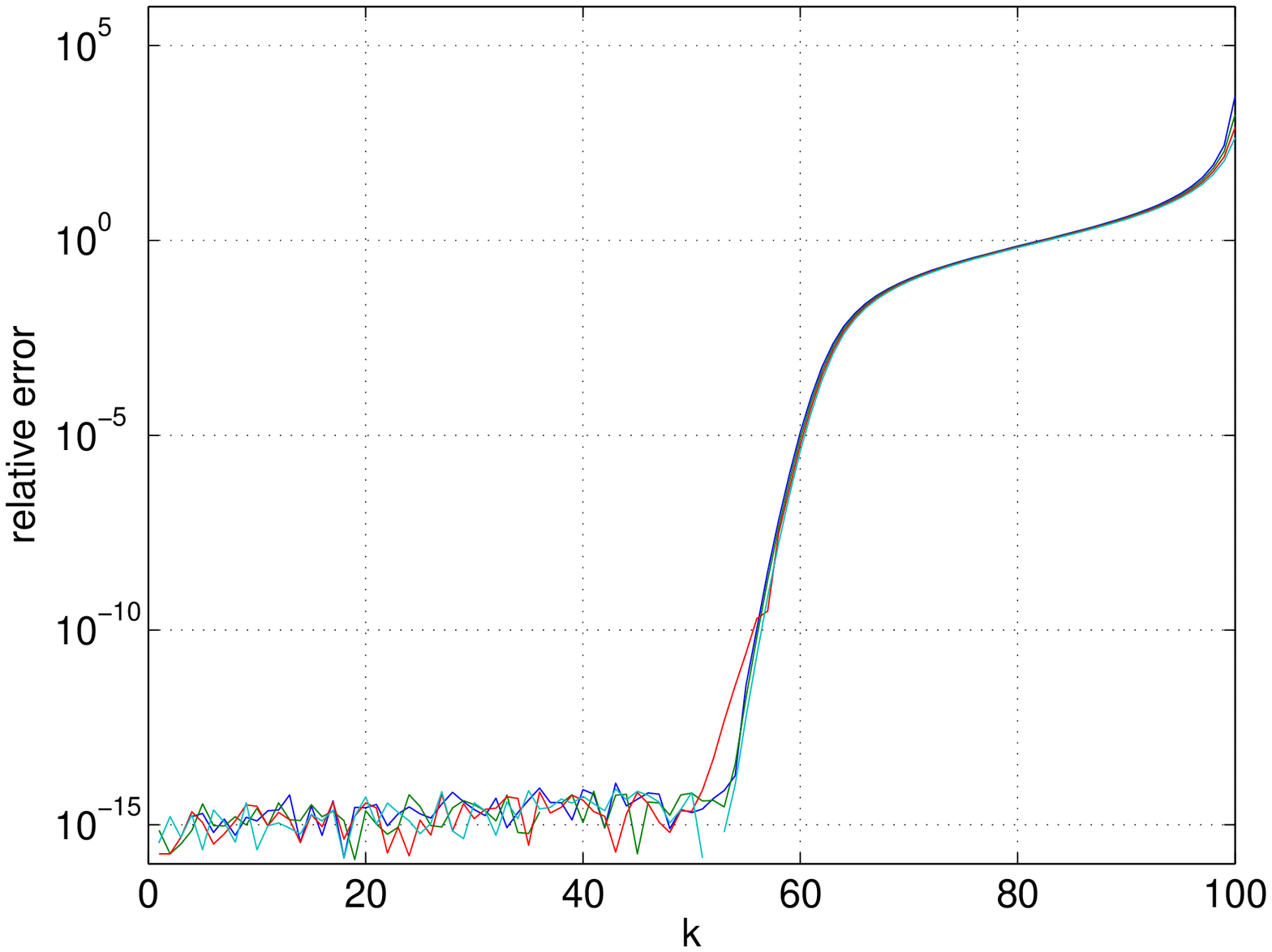}
\includegraphics[width=0.5\textwidth]{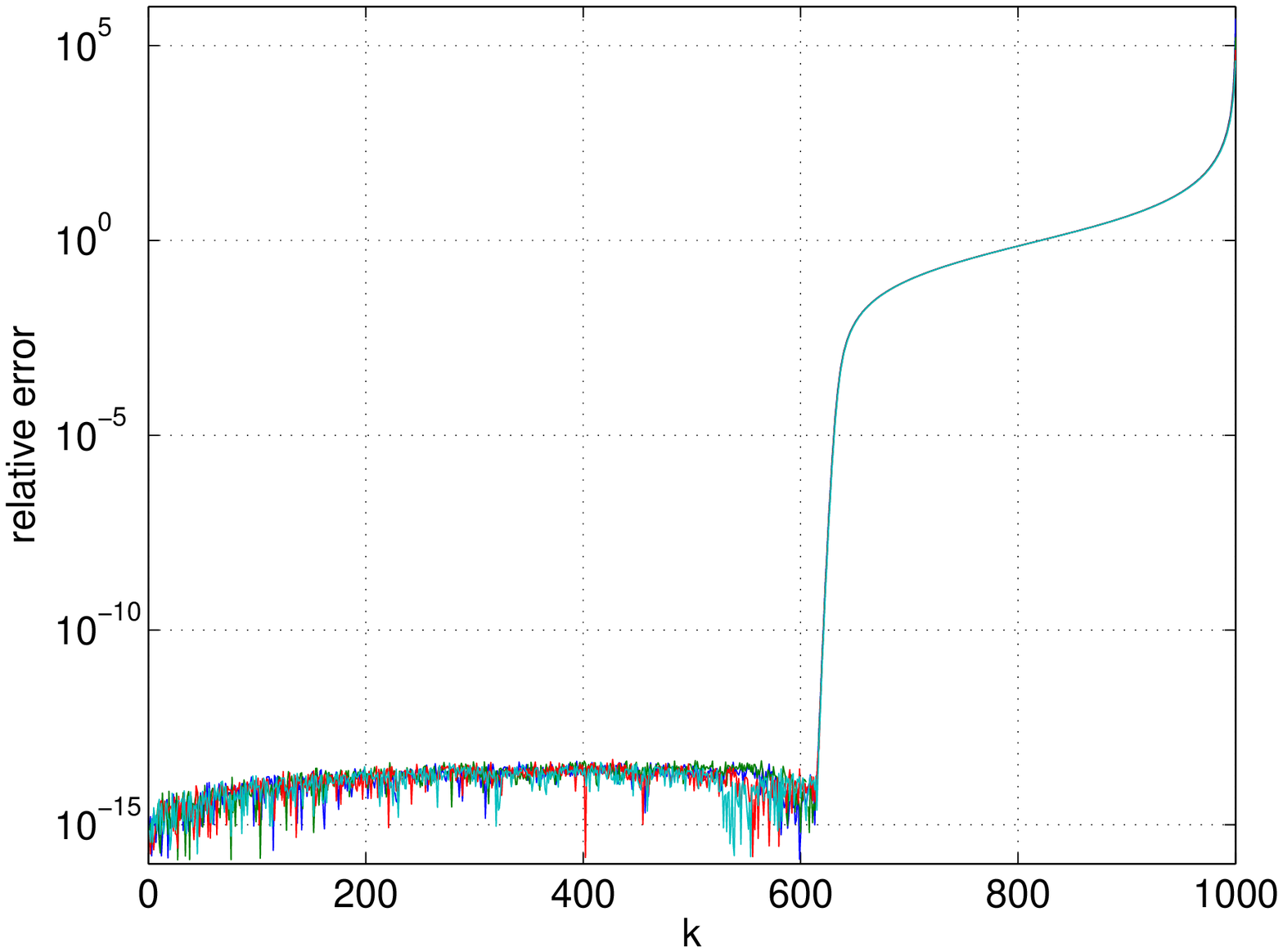}
\caption{Relative error $|\lambda-\lambda_{e}|/|\lambda_e| $ for the entire odd mode spectrum for \texttt{MG}=100 (left) and \texttt{MG} =1000 (right). Exact eigenvalues are $\lambda_e= -k^2 \pi^2$, $k=1,\ldots,$\texttt{MG}. The relative errors for $\gamma=0$, $0.5$, $1$ and $1.5$ are shown but essentially indistinguishable at this scale. }
\label{oddeigs}
\end{figure}

The recurrences (\ref{pmrec}) and (\ref{qmrec}) for the characteristic polynomials can be expressed in the matrix form  
\begin{equation}
\mu \, [p_0(\mu), p_1(\mu), \cdots ] = [p_0(\mu), p_1(\mu), \cdots ] \; M 
\label{recM}
\end{equation}
where the semi-infinite matrix $M$ is tridiagonal plus one row. The matrix $M$ is purely tridiagonal if $\gamma=1/2$ or $3/2$. The roots of the $m$-th order polynomial $p_m(\mu)$ are the eigenvalues of the $m$-by-$m$ matrix $M(0:m-1,0:m-1)$.
A matlab code, \texttt{buildGI2.m}, is provided in appendix \ref{code} which constructs the matrix \texttt{GI2}=$M(0:m,0:m-1)$ for both the even and odd modes by direct implementation of formulas (\ref{pmrec}) and (\ref{qmrec}) with (\ref{Knrec}).  This approach provides an effective and well-conditioned technique to compute the Gegenbauer-Tau eigenvalues as  illustrated in figure \ref{oddeigs} which shows the odd mode eigenvalues for two values of $m\equiv$ \texttt{MG}, the total number of modes, for  $\gamma=0$, 0.5, 1, 1.5, corresponding to Chebyshev-Tau, Legendre-Tau, Chebyshev-Galerkin and Legendre-Tau, respectively. The Gegenbauer-Tau approximations involve an expansion of the solution $u_n(x)$ into odd polynomials up to degree $n=2 \,\tt{MG}+1$ and therefore up to degree 2001 for \texttt{MG}=1000. That calculation shows that slightly more than 60\% of the spectrum is captured with close to machine precision (here double precision IEEE arithmetic), demonstrating the excellent numerical conditioning of the formulation. Comparing the \texttt{MG}=100 and \texttt{MG}=1000 calculations shows that there is a slight ballooning of the round-off error at the higher truncation level. This might be explained by assuming randomness of  roundoff errors with a standard deviation growing like $\sqrt{\texttt{MG}}$. It would be interesting to obtain asymptotic estimates for the high frequency modes, $k > 0.6$ \texttt{MG}. The largest eigenvalues for \texttt{MG}=1000 and $\gamma=0$, 0.5, 1, and 1.5 are, respectively, 
$4.86 \times 10^{12}$, $1.63 \times 10^{12}$, $ 7.61\times 10^{11}$ and $4.07 \times 10^{11}$. So the Legendre-Galerkin method can be said to be slightly less stiff than the other methods. These values are consistent with estimates that the largest eigenvalues are $O($\texttt{MG}$^4)$ \cite{WT88}.

\subsection{Don't Differentiate, Integrate} \label{numerical}

Our approach so far has been theoretical and focused on the basic eigenproblem (\ref{SOP}). For more general two-point boundary value problems, \eg nonlinear problems, it would not be possible to obtain explicit forms such as (\ref{u_n}) for the discrete solution, and the residuals would not be as simple as (\ref{JTR}) or (\ref{RnGeg}). For more general applications it is necessary to select explicit bases for the trial and test functions and to perfom the integrals (\ref{jactau}) or (\ref{jacgal})  by Gauss integration.  

One classical implementation of the Gegenbauer-Tau method is to expand $u_n(x)$ in terms of  Gegenbauer polynomials 
$ %\begin{equation}
u_n(x)= \sum_{l=0}^n a_l \, G_l^{(\gamma)} (x)
$ %\end{equation}
and to use the $n-1$ Gegenbauer polynomials $G_k^{(\gamma)} (x)$, $k=0,\ldots, n-2$ as the test functions in lieu of $p_{n-2}(x)$ in (\ref{jactau}). Those $n-1$ integrals, computed by Gauss quadrature in practice, and the two boundary conditions  yield the $n+1$ equations to determine 
 the $n+1$ coefficients $a_k$.  For the even modes of the simple eigenproblem (\ref{SOP}),  this formulation consists of the even expansion 
 \begin{equation}
 u_{2m}(x)= \sum_{l=0}^m a_l  G_{2l}^{(\gamma)}(x)
 \quad \Rightarrow \quad
 D^2 u_{2m}(x)= \sum_{l=0}^m a_l  \, D^2 G_{2l}^{(\gamma)}(x)
 \label{D2imp}
 \end{equation}
with  the weighted residual equations
 \begin{equation}
\int_{-1}^{1} \left(D^2 u_{2m}- \lambda u_{2m} \right) G_{2k}^{(\gamma)}\, W^{(\gamma)} dx = 0, 
\quad k=0,\ldots, m-1,
\label{wrcg}
 \end{equation}
 where $W^{(\gamma)}=(1-x^2)^{\gamma-1/2}$. Equations (\ref{wrcg}) yield a matrix problem $Aa = \lambda B a$ for $a=[a_0,\ldots,a_m]^T$, where the $m$-by-$(m+1)$ matrix $B$ is diagonal plus one zero column, from orthogonality of the Gegenbauer polynomials (\ref{gegortho}), 
 and the   $m$-by-$(m+1)$ matrix $A$ is upper triangular with $A(k,l)=0$ for $k \ge l$ since, from (\ref{dgegrec2}),  $D^2 G_{2l}(x)$ can be expressed in terms of all the Gegenbauer polynomials of even degree  less than $2l$. The boundary condition $u_n(1)= \sum_{l=0}^m a_l G_{2l}^{(\gamma)}(1) = 0$ allows the elimination of one of the coefficients, $a_0$ or $a_m$ say. This elimination can be expressed in the form 
 $a = C \tilde{a}$ where $\tilde{a}$ is the column vector containing the remaining $m$ coefficients and 
 $C$ is an $(m+1)$-by-$m$ matrix consisting of the $m$-by-$m$ identity matrix plus one full row. This yields the generalized eigenvalue problem $AC \tilde{a} = \lambda BC \tilde{a}$. The structure of the resulting matrices $(AC)$ and $(BC)$ depends on which coefficient is eliminated. If $a_m$ is eliminated, then $(AC)$ is full and $(BC)$ is diagonal. If $a_0$ is eliminated then $(AC)$ is upper triangular and  $(BC)$ is zero everywhere except on the first row and the first lower diagonal. 

\begin{figure}[h]
\begin{center}
\includegraphics[width=0.6\textwidth]{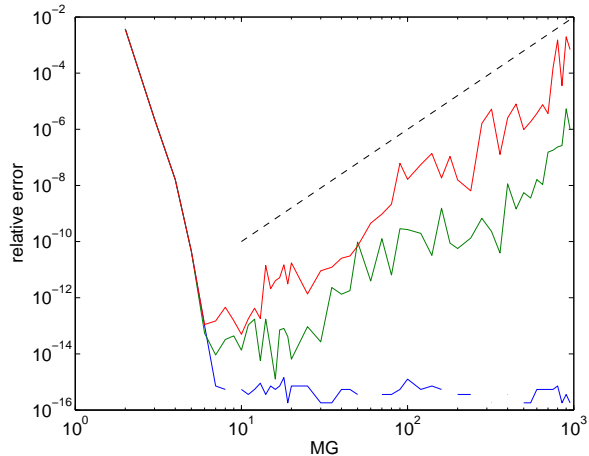}
\caption{Relative error $|\lambda-\lambda_{e}|/|\lambda_e| $ for the first even eigenvalue as a function of \texttt{MG} $\equiv m$ for the Chebyshev-Tau method, $\gamma=0$. The exact eigenvalue $\lambda_e= - \pi^2/4$. The dashed line indicates $MG^4$ scaling of roundoff errors. Three implementations are shown, the differentiation approach (\ref{D2imp}) with $a_m$ eliminated (top curve) and with $a_0$ eliminated (middle curve), and the integration approach (\ref{I2un}). The latter is well-conditioned with errors staying at the level of machine precision $10^{-15}$. The gaps in that curve occur where the approximate eigenvalue is indistinguishable from the numerical value for $\pi^2/4$. }
\end{center}
\label{I2D2comp}
\end{figure}

Many other implementations are possible. For instance, one can use a polynomial expansion that satisfies the boundary conditions \textit{a priori},
$u_{2m}(x)=$ $(1-x^2) \sum_{l=0}^{m-1}  b_l \, \varphi_{2l}(x)$, where $\varphi_{2l}(x)$ is an even polynomial of degree $2l$.  Picking $ \varphi_{2l}(x)= G_{2l}^{(\gamma)}(x)$, the equations (\ref{wrcg}) lead to a generalized eigenvalue problem $A b = \lambda B b$ where this $m$-by-$m$ matrix $A$ is upper triangular and $B$ is tridiagonal.  
All of these formulations are mathematically equivalent; in exact arithmetic they would provide the same eigenvalues as the matrix $M$ in (\ref{recM}). However, the formulations just mentioned use 2nd derivatives of Gegenbauer polynomials and these methods are plagued by roundoff errors that grow like $m^4$, the fourth power of the number of coefficients as illustrated in figure \ref{I2D2comp} \cite{Greengard,TT87}.

There is one formulation that is numerically stable and leads exactly to the tridiagonal plus one row matrix of eqn.\ (\ref{recM}). That formulation consists in expanding not $u_{2m}(x)$ but its 2nd derivative $D^2 u_{2m}(x)$  in terms of Gegenbauer polynomials: 
\begin{equation}
D^2 u_{2m}(x) = \sum_{l=0}^{m-1} c_l \, G_{2l}^{(\gamma)}(x),
%\end{equation} then \begin{equation}
\quad \Rightarrow  \quad 
u_{2m}(x) = \sum_{l=0}^{m-1} c_l \;\SI^2 G_{2l}^{(\gamma)}(x)  + \alpha + \beta  x, 
\label{I2un}
\end{equation}
where $\SI^2$ denotes double integration. That double integration is easily expressed in terms of Gegenbauer polynomials by double integration of the recurrence formulas (\ref{dgegrec2}) which gives
\begin{equation}
\left\{ \begin{aligned}
\SI^2 G_{0}^{(\gamma)}= & \frac{G_{2}^{(\gamma)}}{2(\gamma+1)} + \alpha_0 + \beta_0 x,  \quad
\SI^2 G_{1}^{(\gamma)}= 
\frac{G_{3}^{(\gamma)}}{4(\gamma+1)(\gamma+2)} + \alpha_1 + \beta_1 x,  \\
\SI^2 G_{2}^{(\gamma)}=&  \frac{G_{4}^{(\gamma)}}{4(\gamma+2)(\gamma+3)}
       -\frac{G_{2}^{(\gamma)}}{2(\gamma+1)(\gamma+3)} + \alpha_2 + \beta_2 x, \\ 
\SI^2 G_n^{(\gamma)}=&
\frac{G_{n+2}^{(\gamma)}}{4(\gamma+n+1)(\gamma+ n)} -
\frac{G_{n}^{(\gamma)}}{2(\gamma+ n+1)(\gamma+n-1)} +
\frac{G_{n-2}^{(\gamma)}}{4(\gamma+n)(\gamma+n-1)} \\
& + \alpha_n + \beta_n x.
\end{aligned}
\right.
\label{Igegrec2}
\end{equation}
The constants of integration $\alpha_n$ and $\beta_n$ can be defined arbitrarily since the $\alpha + \beta x$ terms have been included in (\ref{I2un}), so let $\alpha_n=\beta_n=0$ for all $n$. For the even mode expansion considered in this section, we have $\beta =0$ in (\ref{I2un}), so only $\alpha$ survives as the lone constant of integration. That constant is determined from the boundary condition $u_n(1)=0$, which for (\ref{I2un}) reads $ \sum_{l=0}^{m-1} c_l \, \SI^2 G_{2l}^{(\gamma)}(1)  + \alpha =0$. From (\ref{Igegrec2}) with $\alpha_n=\beta_n=0$, one finds that 
\begin{equation}
\alpha = - \sum_{l=0}^{m-1} c_l \, K_{2l}^{(\gamma)}
\end{equation}
with the constants $K_{2l}^{(\gamma)}$ as in (\ref{Kna}) and (\ref{K0K2}). Substituting (\ref{I2un}) with (\ref{Igegrec2}) into (\ref{wrcg}) and using orthogonality of the Gegenbauer polynomials (\ref{gegortho}) yields an eigenvalue problem $Ac = \lambda B c$ where the $m$-by-$m$ matrix $B$ is tridiagonal plus one top row and the $m$-by-$m$ matrix $A$ is diagonal with $A(k,k)= \int_{-1}^1 (G_{2k}^{(\gamma)} )^2 W^{(\gamma)} dx > 0$. The system can thus be rescaled to the form 
\begin{equation}
c = \lambda M c
\label{MGI2}
\end{equation}
 where the matrix $M=A^{-1}B$ is the tridiagonal plus one top row matrix in (\ref{recM}) that was obtained from the characteristic polynomial recurrences (recall that $\mu=1/\lambda$ and that $B$ is tridiagonal if $\gamma=1/2$ or 3/2).  That matrix which consists of the coefficients in (\ref{pmrec}) or (\ref{Igegrec2}) (with $\alpha_n=\beta_n=0$) together with the constants 
 $-K _{2l}^{(\gamma)}$ that modify the first row and impose the boundary condition can now be interpreted as the chopped double Gegenbauer Integration operator with Dirichlet boundary conditions. That is if $f(x) = \sum_{l=0}^{m-1} f_l \, G_{2l}^{(\gamma)} (x)$ then  $g=M^{+} f$ where $M^{+}=M(0:m,0:m-1)$ and  $f=[f_0,f_1,\ldots,f_{m-1}]^T$ provides the $m+1$ even Gegenbauer coefficients of the double integral of  $f(x)$ that vanishes at $x=\pm 1$.  Note that (\ref{pmrec}) and (\ref{MGI2}) provide direct interpretations for the left and right eigenvectors of $M$, respectively.  The problem for the odd modes is entirely analogous and does not need to be repeated here since all the details are available in (\ref{qmrec}) and in the matlab code in appendix \ref{code} which provides \texttt{GI2}= $M^+$. The numerical performance of  two differentiation approaches based on (\ref{D2imp}), and of the integration approach (\ref{I2un}) equivalent to (\ref{recM}), are shown in figure \ref{I2D2comp}  which displays the relative error for the first even mode eigenvalue as a function of $m\equiv$ \texttt{MG} for $\gamma=0$.
 The integration formulation (\ref{I2un}) was proposed by Greengard \cite[p.\ 1077]{Greengard} precisely for the purpose of controlling roundoff errors. This procedure is essentially equivalent to the commonly used reformulation suggested in \cite[p.\ 120]{GO77}, \cite[\S 5.1.2]{CHQZ}.

The Legendre Galerkin (\ie Gegenbauer Tau with $\gamma=3/2$) integration implementation corresponds to Ierley's expansion in associated Legendre polynomials \cite{ierley97}. For (\ref{SOP}), and restricting to even modes, Ierley's expansion consists of $u_n(x)= \sum_{l=0}^{m-1} g_l \, (1-x^2) G_{2l}^{(3/2)}(x)$, where $G_n^{(3/2)}(x) \propto P_n^{(1,1)}(x) \propto D P_{n+1}(x)$ and $P_n(x)$ is the Legendre polynomial of degree $n$ (appendix \ref{orthopolys}). Now  the derivative of eqn.\ (\ref{jaceq}) for $\alpha=\beta=0$ gives $D^2 \left( (1-x^2) D P_{n+1}\right)  = (n+1)(n+2) D P_{n+1}$, and, since $DP_{n+1} \propto G_n^{(3/2)}(x)$,  Ierley' s expansion satisfies
 $D^2 u_n(x) = \sum_{l=0}^{m-1} g_l \, (2l+1)(2l+2) G_{2l}^{(3/2)}(x)$ and corresponds to an expansion of the 2nd derivative of $u_n(x)$ in terms of Gegenbauer polynomials of index $\gamma=3/2$, a special case of the integration approach (\ref{I2un}). Ierley's test functions $(1-x^2) G_{2k}^{(3/2)}(x)$ vanish at $x=\pm 1$, so his equations are (\ref{jacgal}) with $\alpha=\beta=0$ corresponding indeed to a Legendre Galerkin approach (or Gegenbauer Tau with $\gamma=3/2$). This yields an eigenvalue problem of the form $Ag = \lambda B g$ where $A$ is diagonal and $B$ is tridiagonal, where the coefficients $g$ have been renormalized so that $B$ is also symmetric. 

\section{Conclusions}
It has been shown that the eigenvalues of the Jacobi Tau method for
the second derivative operator with Dirichlet boundary conditions
are real, negative and distinct for ranges of the Jacobi indices $\alpha$ and $\beta$. These ranges  include Tau methods with Chebyshev and Legendre polynomials of the 1st and 2nd kinds.
Chebyshev and Legendre Galerkin formulations are included as well but collocation methods are not. Although our work owes much to earlier work by Gottlieb and Lustman \cite{G81,GL83}, we have raised doubts about the validity of their proof for the Chebyshev collocation operator.

Special emphasis has been placed on the symmetric case of the Gegenbauer Tau method where the
range of parameters included in the theorems can be extended and
characteristic polynomials given by successive order approximations 
interlace. The interlacing is between $q_{m-1}(\mu)$ and $p_m(\mu)$ in (\ref{pmqmdefs}), and between 
$p_m(\mu)$ and $q_m(\mu)$, not between $p_m(\mu)$ and $p_{m+1}(\mu)$ or between $q_m(\mu)$ and $q_{m+1}(\mu)$, although we believe the latter hold as well \cite[conjecture 3]{CCW}. Proving such interlacings could allow a proof for the spectrum of the Gegenbauer collocation operator since the parity-reduced residual in that case reads $x DG_n^{(\gamma)}(x)$ which can be written as a linear combination of $G_n^{(\gamma+1)}(x)$ and $G_{n-2}^{(\gamma+1)}(x)$, from (\ref{gegder}) and (\ref{gegrec}).

The characteristic polynomials for Gegenbauer Tau approximations have been shown
to satisfy three term recurrences \emph{plus} a constant term that vanishes for the case of the Legendre Tau and Galerkin methods. Hence for those two  particular cases the
characteristic polynomials are orthogonal, and their roots interlace. A well conditioned matlab code that
computes the roots of the characteristic polynomials for general
Gegenbauer parameter $\gamma$ is  provided in appendix
\ref{code}. In section \ref{numerical}, several 
mathematically equivalent numerical formulations are discussed. The theoretical and practical superiority of the integration method, which is numerically stable, is emphasized. In a forthcoming paper we apply similar methods to the simplified Stokes eigenvalue problem $D^4 u = \lambda D^2 u$ with $u(\pm 1)= Du(\pm 1)=0$ and rigorously identify  classes of spectral methods that are free of spurious eigenvalues. 

\section*{Acknowledgments}

The authors thank Jue Wang for several  helpful calculations in the early stages of this work. 

\appendix

\section{Proof of Theorem \ref{thmP01} and Theorem \ref{thmP02}}
\label{2proofs}

\begin{proof}{(Theorem \ref{thmP01})}
For $A=0$ the theorem reduces to theorem \ref{CCWthm}. Fix now $A>0$
but otherwise arbitrary. Let
\begin{equation}
f_n(x;\mu)=
\sum_{k=0}^{n}\mu^kD^{k}P_n^{(\alpha,\beta)}(x)
+A\sum_{k=0}^{n-1}\mu^kD^{k}P_{n-1}^{(\alpha,\beta)}(x)  \end{equation}
with $\mu$ such that $f_n(1;\mu)=0$ and $D:=d/dx$. Then
$f_n(x;\mu)$ satisfies the following differential equation
\begin{equation}
\left(f_n-P_{n}^{(\alpha,\beta)}-AP_{n-1}^{(\alpha,\beta)}\right)=\mu
\frac{d{f_n}}{d{x}}.
\end{equation}
Multiplying by $\frac{d f_n^*(x,\mu)}{d x} (1+x)$, integrating from
$-1$ to $1$ in the Jacobi norm and adding the conjugate we obtain:
\begin{multline}   \label{equation1}
%\begin{gathered}
\int_{-1}^1\frac{d{|f_n|^2}}{d{x}}(1+x) \jacw8 dx -
\int_{-1}^1\left(\frac{d{f_n}}{d{x}}+\frac{d{f^*_n}}{d{x}}\right)(1+x)P_{n}^{(\alpha,\beta)}(x)
\jacw8 dx    \\
-A\int_{-1}^1\left(\frac{d{f_n}}{d{x}}+\frac{d{f^*_n}}{d{x}}\right)(1+x)P_{n-1}^{(\alpha,\beta)}(x)
\jacw8 dx = (\mu+\mu^*)\int_{-1}^1
\left|\frac{d{f_n}}{d{x}}\right|^2(1+x) \jacw8 dx.
%\end{gathered} \end{equation}
\end{multline}
For the first term of equation (\ref{equation1}),  integration by parts yields
\begin{equation}
\int_{-1}^1\frac{d{|f_n|^2}}{d{x}}(1+x) \jacw8 dx= -\int_{-1}^1
|f_n|^2 \frac{\jacw8}{(1-x)} \left(\beta+1-\alpha-(\beta+1+\alpha)x \right)dx.
\end{equation}
For $\beta>-1$ and $\alpha \le 0$ the factor
$\beta+1-\alpha-(\beta+1+\alpha)x$ is nonnegative for all $x \in
[-1,1]$.  The second term of (\ref{equation1}) can be expanded as 
%\begin{equation} \begin{gathered}
\begin{multline}
\int_{-1}^1\left(\frac{d{f_n}}{d{x}}+\frac{d{f^*_n}}{d{x}}\right)(1+x)P_{n}^{(\alpha,\beta)}
\jacw8 dx=2\int_{-1}^1x DP_{n}^{(\alpha,\beta)} P_{n}^{(\alpha,\beta)} \jacw8dx=\\
2B_{n-1}\int_{-1}^1x
P_{n-1}^{(\alpha,\beta)}P_{n}^{(\alpha,\beta)}\jacw8
dx=\frac{2B_{n-1}a_{1,n-1}}{a_{3,n-1}}\int_{-1}^1
\left(P_{n}^{(\alpha,\beta)}\right)^2 \jacw8
dx=\frac{2B_{n-1}a_{1,n-1}}{a_{3,n-1}} h_n^{\alpha,\beta},
%\end{gathered}\end{equation}
\end{multline}
where we have used expression (\ref{Jacobider}),  the Jacobi
recurrence relation (\ref{jacrec}) and orthogonality of 
$\jacx$ to all polynomials of degree less than $n$
with respect to the Jacobi weight $\jacw8(x)$. Similarly, the third term of equation (\ref{equation1}) can be calculated as 
%\begin{equation} \begin{gathered}
\begin{multline}
\int_{-1}^1\left(\frac{d{f_n}}{d{x}}+\frac{d{f^*_n}}{d{x}}\right)(1+x)P_{n-1}^{(\alpha,\beta)}\jacw8dx  \\
 = 2\int_{-1}^1 (1+x) DP_{n}^{(\alpha,\beta)}
P_{n-1}^{(\alpha,\beta)} \jacw8 dx +(\mu+\mu^*) \int_{-1}^1 x D^2
P_{n}^{(\alpha,\beta)} P_{n-1}^{(\alpha,\beta)} \jacw8dx   \hfill   \\ 
+ \, 2A \int_{-1}^1 x DP_{n-1}^{(\alpha,\beta)} P_{n-1}^{(\alpha,\beta)} \jacw8dx   \\
 = 2B_{n-1}\int_{-1}^1 \left(P_{n-1}^{(\alpha,\beta)}\right)^2 \jacw8
dx + 2B_{n-1}\int_{-1}^1 x \left(P_{n-1}^{(\alpha,\beta)}\right)^2 \jacw8 dx  \hfill  \\ 
\hfill  + (\mu+\mu^*)B_{n-1}B_{n-2}\int_{-1}^1 x
P_{N-2}^{(\alpha,\beta)}P_{n-1}^{(\alpha,\beta)} \jacw8dx
+2A B_{n-2} \int_{-1}^1 xP_{N-2}^{(\alpha,\beta)}P_{n-1}^{(\alpha,\beta)} \jacw8 dx  \\
= 2B_{n-1} h_{n-1}^{\alpha,\beta}-2\frac{a_{2,n-1}}{a_{3,n-1}} B_{n-1} h_{n-1}^{\alpha,\beta}
+(\mu+\mu^*) \frac{a_{1,n-2}}{a_{3,n-2}} B_{n-1}B_{n-2}
h_{n-1}^{\alpha,\beta}+2A\frac{a_{1,n-2}}{a_{3,n-2}}B_{n-2}h_{n-1}^{\alpha,\beta}. \hfill
\end{multline}
%\end{gathered} \end{equation}
%
Substituting these expressions back into equation (\ref{equation1}) yields
\begin{multline}
-\left[\int_{-1}^1 |f_n|^2
\frac{\jacw8}{(1-x)}(\beta+1-\alpha-(\beta+1+\alpha)x)dx+\frac{2B_{n-1}a_{1,n-1}}{a_{3,n-1}}
h_n^{\alpha,\beta}  \;  \qquad   \right. \\
  \left. 
+ \; 2A\left( \left(1-\frac{a_{2,n-1}}{a_{3,n-1}}\right) B_{n-1}
+A\frac{a_{1,n-2}}{a_{3,n-2}}B_{n-2} \right)
h_{n-1}^{\alpha,\beta}\right]  \qquad  \\ 
\hfill  = (\mu+\mu^*) \left[\int_{-1}^1
\left|\frac{d{f_n}}{d{x}}\right|^2(1+x) \jacw8
dx+A\frac{a_{1,n-2}}{a_{3,n-2}} B_{n-1}B_{n-2}
h_{n-1}^{\alpha,\beta}\right].  
\end{multline}
Since $1>\frac{a_{2,n-1}}{a_{3,n-1}}$ the left-hand side is positive
and $\Re(\mu)<0$, which ensures stability.
\end{proof}

\bigskip

\begin{proof}{(Theorem \ref{thmP02})}
For $A=0$ the theorem reduces again to theorem \ref{CCWthm}. Fix now
$A>0$. Let
\begin{equation}\label{fN2b}
f_n(x;\mu)=\sum_{k=0}^{N}\mu^kD^{k}P_n^{(\alpha,\beta)}(x)+A\mu^2
\sum_{k=0}^{N-1}\mu^kD^{k}P_{n-1}^{(\alpha,\beta)}(x).
\end{equation}
with $f_n(1;\mu)=0$. Then $f_n(x;\mu)$ satisfies the differential equation
\begin{equation}
\frac{1}{\mu}\left(f_n-P_{n}^{(\alpha,\beta)}-A\mu^2
P_{n-1}^{(\alpha,\beta)}\right)= \frac{d{f_n}}{d{x}}.
\end{equation}
Multiplying by $f_n^*(x,\mu) \frac{(1+x)}{(1-x)}$, integrating from
$-1$ to $1$ and adding the conjugate yields
%\begin{equation} \begin{gathered}
\begin{multline}
 \label{equation1d}
\left(\frac{1}{\mu}+\frac{1}{\mu^*}\right)\int_{-1}^1|f_n|^2
\frac{(1+x)}{(1-x)} \jacw8 dx - \frac{1}{\mu}\int_{-1}^1 f^*_n
\frac{(1+x)}{(1-x)} P_{n}^{(\alpha,\beta)}(x) \jacw8 dx   \hfill  \\ 
\qquad \qquad \qquad - \frac{1}{\mu^*}\int_{-1}^1 f_n \frac{(1+x)}{(1-x)}
P_{n}^{(\alpha,\beta)}(x) \jacw8 dx -A\mu \int_{-1}^1 f^*_n
\frac{(1+x)}{(1-x)}P_{n-1}^{(\alpha,\beta)}(x) \jacw8 dx \qquad
 \\
\hfill  -A\mu^* \int_{-1}^1 f_n
\frac{(1+x)}{(1-x)}P_{n-1}^{(\alpha,\beta)}(x) \jacw8 dx
\; =\int_{-1}^1 \frac{d{|f_n|^2}}{d{x}} \frac{(1+x)}{(1-x)} \jacw8 dx.
\end{multline}
%\end{gathered}\end{equation}
Integration by parts on the first term gives
\begin{equation}
\int_{-1}^1\frac{d{|f_n|^2}}{d{x}}(1+x) \jacw8 dx= -\int_{-1}^1
|f_n|^2 \frac{\jacw8}{(1-x)^2} \left(\beta-\alpha+2-(\beta+\alpha)x\right)dx.
\end{equation}
For $\beta>-1$ and $\alpha \le 1$ the factor
$\beta-\alpha+2-(\beta+\alpha)x$ is nonnegative for all $x \in
[-1,1]$.  For the other terms on the left hand side of (\ref{equation1d}), recall that $f_n(1;\mu)=0$ so write
\begin{equation}\label{eqnfN2b}
f_n=(1-x)\sum_{k=0}^{n-1} c_k P_{k}^{(\alpha,\beta)}(x)
\end{equation}
then
\begin{multline}
\int_{-1}^1  \frac{(1+x)}{(1-x)}f^*_n P_{n}^{(\alpha,\beta)}(x)
\jacw8 dx=\int_{-1}^1 c_{n-1}^* x
P_{n-1}^{(\alpha,\beta)}(x)P_{n}^{(\alpha,\beta)}(x) \jacw8 dx
\\
= c_{n-1}^*\frac{a_{1,n-1}}{a_{3,n-1}} \int_{-1}^1
\left(P_{n}^{(\alpha,\beta)}(x)\right)^2 \jacw8
dx=c_{n-1}^*\frac{a_{1,n-1}}{a_{3,n-1}} h_n^{\alpha,\beta}.
\end{multline}
Also
%\begin{equation} \begin{gathered}
\begin{multline}
\int_{-1}^1  \frac{(1+x)}{(1-x)}f^*_n P_{n-1}^{(\alpha,\beta)}(x) \jacw8 dx   \\ 
 =\int_{-1}^1 c_{n-1}^*(1+x) P_{n-1}^{(\alpha,\beta)}(x)
P_{n-1}^{(\alpha,\beta)}(x) \jacw8 dx +\int_{-1}^1 c_{n-2}^* x
P_{n-2}^{(\alpha,\beta)}(x) P_{n-1}^{(\alpha,\beta)}(x) \jacw8  dx \hfill  \\
= c_{n-1}^* \int_{-1}^1 \left(P_{n-1}^{(\alpha,\beta)}(x)\right)^2
\jacw8 dx-c_{n-1}^* \frac{a_{2,n-1}}{a_{3,n-1}}\int_{-1}^1
\left(P_{n-1}^{(\alpha,\beta)}(x)\right)^2 \jacw8 dx  \hfill \\
\hfil + c_{n-2}^*\frac{a_{1,n-2}}{a_{3,n-2}} \int_{-1}^1
\left(P_{n-1}^{(\alpha,\beta)}(x)\right)^2 \jacw8 dx \\
 =\left(c_{n-1}^*
\left[1-\frac{a_{2,n-1}}{a_{3,n-1}}\right]+c_{n-2}^*\frac{a_{1,n-2}}{a_{3,n-2}}
\right) h_{n-1}^{\alpha,\beta}.  \hfill 
\end{multline}
%\end{gathered} \end{equation}
%
%
Explicit values of $c_{n-1}$ and $c_{n-2}$ follow from 
equation (\ref{eqnfN2b}) 
\begin{multline} %\begin{split}
f_n= (1-x)\sum_{k=0}^{n-1} c_k P_{k}^{(\alpha,\beta)}(x) \\
= c_{n-1}
P_{n-1}^{(\alpha,\beta)}(x)-c_{n-1}
xP_{n-1}^{(\alpha,\beta)}(x)-c_{n-2}
xP_{n-2}^{(\alpha,\beta)}(x)+O(n-2) \\
=-c_{n-1}\frac{a_{1,n-1}}{a_{3,n-1}}P_{n}^{(\alpha,\beta)}(x)+\left(c_{n-1}\left[1+\frac{a_{2,n-1}}{a_{3,n-1}}
\right] -c_{n-2}
\frac{a_{1,n-2}}{a_{3,n-2}}\right)P_{n-1}^{(\alpha,\beta)}(x)+O(n-2).
\end{multline}
%\end{split} \end{equation}
%
Now from equation (\ref{fN2b}) 
\begin{equation}
\begin{split}
f_n=& P_{n}^{(\alpha,\beta)}(x)+\mu DP_{n}^{(\alpha,\beta)}(x)+ A\mu^2
P_{n-1}^{(\alpha,\beta)}(x)+O(n-2) \\
=& P_{n}^{(\alpha,\beta)}(x)+B_{n-1} \mu
P_{n-1}^{(\alpha,\beta)}(x) + A \mu^2 P_{n-1}^{(\alpha,\beta)}(x)
+O(n-2).
\end{split}
\end{equation}
Comparing these two expressions for $f_n$ gives
\begin{equation}
c_{n-1}=-\frac{a_{3,n-1}}{a_{1,n-1}}  \qquad \qquad
c_{n-2}=\frac{a_{3,n-2}}{a_{1,n-2}}
\left(-\frac{a_{2,n-1}}{a_{1,n-1}}-\frac{a_{3,n-1}}{a_{1,n-1}}-B_{n-1}\mu-A
\mu^2 \right)
\end{equation}
Substituting all these results back into  (\ref{equation1d}) yields
%
%\begin{equation} \begin{gathered}
\begin{multline}
\left(\frac{1}{\mu}+\frac{1}{\mu^*}\right)\left[\int_{-1}^1|f_n|^2
\frac{(1+x)}{(1-x)} \jacw8 dx+h_n^{\alpha,\beta} \right] +
2A(\mu+\mu^*) \frac{a_{3,n-1}}{a_{1,n-1}} h_{n-1}^{\alpha,\beta}  \hfill \\
\hfill \quad + 2AB_{n-1} |\mu|^2 h_{n-1}^{\alpha,\beta}+(\mu+\mu^*) A^2 |\mu|^2 h_{n-1}^{\alpha,\beta}
=-\int_{-1}^1|f_n|^2
\frac{(\alpha-\beta+2-(\alpha+\beta)x)}{(1-x)^2} \jacw8 dx
\end{multline}
%\end{gathered}\end{equation}
%
or after rearranging some of the terms 
\begin{multline}%\begin{equation}\begin{gathered}
(\mu+\mu^*)\left(\frac{1}{|\mu|^2}\left[\int_{-1}^1|f_n|^2
\frac{(1+x)}{(1-x)} \jacw8 dx+ h_n^{\alpha,\beta}\right] 
+\big(2A \frac{a_{3,n-1}}{a_{1,n-1}}+ A^2 |\mu|^2 \big)h_{n-1}^{\alpha,\beta} \right) \\
=-\int_{-1}^1|f_n|^2
\frac{(\beta-\alpha+2-(\beta+\alpha)x)}{(1-x)^2} \jacw8 dx-
2AB_{n-1} |\mu|^2 h_{n-1}^{\alpha,\beta}. 
\end{multline}%\end{gathered}\end{equation}
The right hand side is negative so this implies that  $\Re(\mu)<0$.
\end{proof}

\section{Jacobi and Gegenbauer polynomials}\label{orthopolys}

\subsection{Jacobi Polynomials}

\label{Jacpolys}

The Jacobi polynomials $\jacx$ are suitably standardized
orthogonal polynomials on the interval $(-1,1)$, with weight
function $\jacw8=(1-x)^{\alpha}(1+x)^{\beta}$. The class of Jacobi
polynomials $\jacx$ includes  Gegenbauer (Ultraspherical) polynomials when $\alpha=\beta$, Chebyshev polynomials when $\alpha=\beta=-1/2$ and Legendre polynomials when $\alpha=\beta=0$. 
\begin{definition}
\label{def:jac} The {\it Jacobi polynomial},
$P_n^{(\alpha,\beta)}(x)$, of degree $n$, can be defined by
\begin{equation}
P_n^{(\alpha,\beta)}(x):=\frac{1}{2^n}\sum_{k=0}^n\ \binom{n+\alpha}{k} \binom{n+\beta}{n-k}(x-1)^{n-k}(x+1)^{k},\quad
\alpha,\,\, \beta >-1,  \label{jacobi}
\end{equation}
\end{definition}
where the binomial coefficient $\binom{\alpha}{k} = (\alpha)
(\alpha-1) \cdots (\alpha-k+1)/k!$.  Jacobi polynomials are the most
general class of polynomial solutions of a singular Sturm-Liouville
problem on the interval $-1<x<1$
% \ie an eigenvalue problem of the form \cite[\S 9.2]{CHQZ}
%\begin{equation} -\left(p(x) \psi' \right)' + q(x) \psi = \lambda w(x) \psi, \label{SL} \end{equation}
%with $\psi'(\pm 1)$ bounded, where $p(x)$, $q(x)$ and $w(x)$ are
%continuously differentiable non-negative functions on the open
%interval $(-1,1)$ with $p(\pm 1)=0$ and $p(x)>0$ if $|x| \neq 1$.
and this is directly related to their excellent approximation
properties \cite[\S 9.2.2, \S 9.6.1]{CHQZ}. %, up to normalization constants.
The Jacobi polynomial $\jacx$  satisfies the differential equation
\begin{equation}
\frac{d}{dx} \left( (1-x)^{\alpha+1} (1+x)^{\beta+1} \frac{d}{dx} y \right) = 
n(n+\alpha+\beta+1) (1-x)^{\alpha}(1+x)^{\beta} y.
\label{jaceq}
\end{equation}
Jacobi polynomials (\ref{jacobi}) are orthogonal with respect to the weight 
$W_{\alpha,\beta}(x)=(1-x)^{\alpha}(1+x)^{\beta}$
 \begin{equation}
\int_{-1}^{1}  (1-x)^{\alpha}(1+x)^{\beta} \,P_m^{(\alpha,\beta)}
P_n^{(\alpha,\beta)} \,dx=\left\{%
\begin{array}{ll}
    0, & m \ne n,\\
    h_n^{\alpha,\beta}, & m = n, 
\end{array}%
\right. \label{jacortho}
\end{equation}
where
\begin{equation}
h_n^{\alpha,\beta}=\frac{2^{\alpha+\beta+1}}{2n+\alpha+\beta+1}
\frac{\Gamma{(n+\alpha+1)} \Gamma{(n+\beta+1})}{n!
\Gamma{(n+\alpha+\beta+1)}}.
\end{equation}
Orthogonal polynomials satisfy a three term recurrence relation, for  the Jacobi polynomials this reads
\begin{multline} %\begin{equation} \begin{split}
\label{jacrec}
2(n+1)(n+\alpha+\beta+1)(2n+\alpha+\beta)
P_{n+1}^{(\alpha,\beta)}(x)= \\
\left((2n+\alpha+\beta+1)(\alpha^2-\beta^2)+(2n+\alpha+\beta)_3\;x
\right)P_{n}^{(\alpha,\beta)}(x) \\
 -2(n+\alpha)(n+\beta)(2n+\alpha+\beta+2)P_{n-1}^{(\alpha,\beta)}(x).
\end{multline} %\end{split}\end{equation}
where
$(2n+\alpha+\beta)_3=(2n+\alpha+\beta)(2n+\alpha+\beta+1)(2n+\alpha+\beta+2)$.
To ease the notation in calculations we write the
recurrence relation in the form
\begin{equation}
a_{1,n}P_{n+1}^{(\alpha,\beta)}(x)=(a_{2,n}+a_{3,n}x)P_{n}^{(\alpha,\beta)}(x)-a_{4,n}P_{n-1}^{(\alpha,\beta)}(x).
\end{equation}
Two other useful relations involving derivatives of Jacobi polynomials \cite{Doha} are
\begin{equation}\label{jacder}
\frac{d}{dx}P_n^{(\alpha,\beta)}(x)=\frac{1}{2}(n+\alpha+\beta+1)P_{n-1}^{(\alpha+1,\beta+1)}(x).
\end{equation}
and
\begin{equation}\label{Jacobider}
\frac{d}{dx} P_{n+1}^{(\alpha,\beta)}(x)=B_n
P_{n}^{(\alpha,\beta)}(x)+p_{n-1}(x)
\end{equation}
with
$B_n=\frac{(2n+\alpha+\beta+1)(2n+\alpha+\beta+1)}{(n+\alpha+\beta+1)}$
and $p_{n-1}(x)$ a polynomial of degree $n-1$.
%
%The Jacobi polynomials with $\alpha = \beta =\gamma-\frac{1}{2}$  correspond to Gegenbauer polynomials \cite[22.5.20]{AS}
%\begin{equation}
%P_n^{(\gamma-\frac{1}{2},\gamma-\frac{1}{2})}(x)=\frac{\Gamma{(2\gamma)}
%\Gamma{(n+\gamma+\frac{1}{2})}}{\Gamma{(n+2\gamma)}\Gamma{(\gamma+\frac{1}{2})}} C_n^{(\gamma)}(x).  \label{jacgeg}
%\end{equation}

\subsection{Gegenbauer Polynomials}

\label{Gegpolys}

The Gegenbauer (a.k.a.\ Ultraspherical) polynomials
$C_n^{(\gamma)}(x)$, $\gamma > -1/2$, of degree $n$ are the Jacobi
polynomials  with $\alpha=\beta=\gamma-1/2$, up to normalization  \cite[22.5.20]{AS}.
They are symmetric (even for $n$ even and odd for $n$ odd) orthogonal polynomials with weight
function $W(x)=(1-x^2)^{\gamma-\frac{1}{2}}$. Since the standard
normalization \cite[22.3.4] {AS}, is singular for the Chebyshev case
$\gamma=0$, we use a non-standard normalization that includes
the Chebyshev case but preserves the simplicity of the Gegenbauer
recurrences. Set
\begin{equation}
G_0^{(\gamma)}(x):=1, \qquad G_n^{(\gamma)}(x) := \frac{C_n^{(\gamma)}(x)}{2\gamma}, \quad  n \ge 1.
\end{equation}
We refer to these non-standard Gegenbauer polynomials as  ns-Gegenbauer for short.
The ns-Gegenbauer polynomials satisfy the orthogonality relationship
\begin{equation}
\int_{-1}^{1}  (1-x^2)^{\gamma-1/2} \,G_m^{(\gamma)}
G_n^{(\gamma)} \,dx =\left\{%
\begin{array}{ll}
    0, & m \ne n, \\
    h_n^{\gamma}, & m = n,
\end{array}%
\right.    \label{gegortho}
\end{equation}
where
 \cite[22.2.3]{AS}, 
\begin{equation}
h_n^{\gamma} = \frac{\pi  2^{-1-2 \gamma} \Gamma(n+2
\gamma)}{\gamma^2 (n+\gamma)n! \Gamma^2(\gamma)}.
\end{equation}
The derivative recurrence formula (\ref{jacder}) for ns-Gegenbauer  polynomials
reads
\begin{equation}
\label{gegder} \frac{d}{dx} G_{n+1}^{(\gamma)} = 2 (\gamma+1) G_n^{(\gamma+1)},
\end{equation}
(for $C_n^{(\gamma)}$ this is formula \cite[A.57]{Boyd}), 
and the three-term recurrence takes the simple form % for ns-Gegenbauer polynomials is
\begin{equation}
(n+1)G_{n+1}^{(\gamma)}= 2(n+\gamma) x G_{n}^{(\gamma)} - (n-1 + 2
\gamma) G_{n-1}^{(\gamma)}, \quad n\ge 2,
 \label{gegrec},
\end{equation}
with % FW 06/07/6
\begin{equation}
G_{0}^{(\gamma)}(x) =1, \quad G_{1}^{(\gamma)}(x) =x, \quad
G_2^{(\gamma)}=(\gamma+1)x^2-\frac{1}{2}.
\end{equation}
Differentiating the recurrence (\ref{gegrec})  with respect to $x$
and  subtracting from the corresponding recurrence for $\gamma+1$
using (\ref{gegder}), yields \cite[22.7.23]{AS} 
\begin{equation} \label{ggp}
(n+\gamma) G_{n}^{(\gamma)} =  (\gamma+1)
 \left[G_{n}^{(\gamma+1)} - G_{n-2}^{(\gamma+1)} \right], \qquad n \ge 3.
\end{equation}
Combined with  (\ref{gegder}), this leads to the important derivative
recurrence between ns-Gegenbauer polynomials of same index $\gamma$
\begin{eqnarray}
G_0^{(\gamma)}(x) = \frac{d}{dx} G_1^{(\gamma)}(x),& \quad
2 (1+\gamma) G_1^{(\gamma)}(x) = \frac{d}{dx} G_2^{(\gamma)}(x), \nonumber \\
2 (n+\gamma) G_{n}^{(\gamma)} =\; &  \frac{d}{dx}
 \left[G_{n+1}^{(\gamma)} - G_{n-1}^{(\gamma)} \right]. \hspace{54pt}
 \label{dgegrec}
\end{eqnarray}
Evaluating the Gegenbauer polynomial at $x=1$ we find \cite[22.4.2]{AS}, 
\begin{equation}\label{Gat1}
G_n^{(\gamma)}(1)=\frac{1}{2\gamma}C_n^{(\gamma)}(1)=\frac{1}{2\gamma}
\binom{2\gamma + n-1}{n}
\end{equation}
where $\binom{2\gamma +  n -1}{n}=
(2\gamma+n-1)(2\gamma+n-2)\cdots (2\gamma)/n!
=\frac{\Gamma{(2\gamma + n)}}{n! \Gamma{(2\gamma)}}$.

Gegenbauer polynomials correspond to Chebyshev polynomials of
the 1st kind, $T_n(x)$, when $\gamma=0$, to Legendre $P_n(x)$ for
$\gamma=1/2$ and to Chebyshev of the 2nd kind, $U_n(x)$, for
$\gamma=1$. For the non standard normalization, 
\begin{equation}\label{limit}
G_n^{(0)}(x) =\frac{T_n(x)}{n} , \quad G_n^{(1/2)}(x)=P_n(x), \quad
G_n^{(1)}(x)=\frac{U_n(x)}{2}.
\end{equation}

\section{Matlab code for Gegenbauer-Tau Double Integration} \label{code}  \begin{verbatim}
function GI2=buildGI2(MG,g,ip)
% buildGI2 produces the Gegenbauer-Tau double integration operator GI2 with
% Dirichlet boundary conditions u(+/-1)=0 for even (ip=0) or odd (ip=1) solutions.
%
%  GI2 = buildGI2(MG,g,ip) yields the (MG+1)-by-MG tridiagonal + 1 row matrix GI2.
%  (2*MG+ip) is the degree of the polynomial expansion, g is the Gegenbauer index
%  g=0 is Chebyshev-Tau, g=1/2 is Legendre-Tau,  g=1 is Chebyshev Galerkin,
%  g=3/2 is Legendre-Galerkin. g must be greater than -1/2.
%
% EXAMPLE:  Cheb-Galerkin odd mode eigenvalues compared to exact values:
%      MG=20; GI2=buildGI2(MG,1,1); M=GI2(1:end-1,:); eCG=sort(1./abs(eig(M)));
%      k=[1:MG]; semilogy(k,k.^2*pi^2,k,eCG,'o')
%
% Fabian Waleffe & Marios Charalambides, 2005, 2006

n=2*(1:MG-1)+ip; 

dm=1./(4*(g+n+1).*(g+n));d0=-1./(2*(g+n+1).*(g+n-1));dp=1./(4*(g+n).*(g+n-1));

T=diag(dm(1:MG-2),-1)+diag(d0)+diag(dp(2:MG-1),1);  % Tridiagonal part

% K_n by recurrence (minus sign included)
if (MG>2), Kn=zeros(1,MG-2); K3=(2*g-1)*(3-2*g)/120;
if     (ip==0), Kn(1)=(4*g^2-1)*(3-2*g)/720; %m=2, n=4, Kn(m)=K_{2m+2}
elseif (ip==1), Kn(1)=K3*(2*g+2)*(2*g+1)/42; %m=2, n=5, Kn(m)=K_{2m+3}
else error('    ip must be 0 or 1'), end 
for m=2:MG-2;
    n=2*m+ip;  Kn(m)=Kn(m-1)*(2*g+n-1)*(2*g+n-2)/((n+4)*(n+3));
end, end

% 1st row and 1st column
if (ip==0), M00=-(2*g+1)/(4*g+4);
   M01=(7-g-2*g^2)*(1+2*g)/(48*(2+g)*(1+g)); M10=1/(2*g+2);
elseif (ip==1), M00=-(2*g+1)/(12*g+24);
   M01=1/(4*(g+3)*(g+2)) + K3;  M10=1/(4*(g+1)*(g+2));
end

r1=[M00, M01, Kn]; c1=[M10; zeros(MG-2,1)]; re=[zeros(1,MG-1),dm(end)];
GI2=[r1; c1,T; re];
\end{verbatim}

\bibliographystyle{siam}
\bibliography{CWJacTau2.bib}
\end{document}